\documentclass[12pt,leqno]{amsart}

\usepackage[utf8]{inputenc}   % LaTeX, comprends les accents !
\usepackage[T1]{fontenc}      % Police contenant les caract??res français
\usepackage{geometry}         % Définir les marges

\usepackage{amssymb}
\usepackage{amsfonts}
\usepackage{amscd}
\usepackage[matrix]{xypic}
\usepackage{enumerate}
\voffset-1cm \markleft{\rm Poisson} 
\newtheorem{theorem}{Theorem}
\newtheorem{definition}[theorem]{Definition}
\newtheorem{corollary}[theorem]{Corollary}

\newtheorem{lemma}[theorem]{Lemma}

\newtheorem{proposition}[theorem]{Proposition}
\newtheorem{remark}[theorem]{Remark}

\newtheorem{defi}{D\'{e}finition}

\newcommand{\ant}{anti-associative }

\newcommand{\K}{\mathbb K}

\newcommand{\N}{\mathbb N}

\newcommand{\KS}{{\mathbb K}[\Sigma_3]}
\newcommand{\ST}{\Sigma_3}

\newcommand{\ra}{\rightarrow}

\newcommand{\pf}{\noindent{\it Proof. }}
\newcommand{\ds}{\displaystyle}

\newcommand{\im}{{\rm Im }}

\newcommand{\bu}{\bullet}

\newcommand\rk{\text{\rm rk}}

\newcommand{\Id}{\mathrm{id}}

\newcommand{\LE}{\mathcal{L}}
\newcommand{\vv}{a_1Id+a_2\tau_{12}+a_3\tau_{13}+a_4\tau_{23}+a_5c+a_6c^2}

\usepackage{graphicx}

\title{Deformation quantization of nonassociative algebras}
\author{ Elisabeth Remm $^{1,2}$}
\date{}
\address{$^1$ Universit\'e de Haute Alsace. IRIMAS UR 7499 - Department of Mathematics, F-68100 Mulhouse, France \quad \quad 
$^2$ Université de Strasbourg}
\email{elisabeth.remm@uha.fr}

\begin{document}

\maketitle

\noindent{\bf Abstract.} We investigate formal deformations of certain classes of nonassociative algebras including classes of $\KS$-associative algebras,  Lie-admissible algebras and anti-associative algebras.    In a process which is similar to Poisson algebra for the associative case we identify  for each type of algebras $(A, \mu)$ a type of algebras   $(A,\mu, \psi)$ such that  formal deformations of $(A,\mu)$ appear  as quantizations  of $(A,\mu, \psi)$. The process of polarization/depolarization associate to each  nonassociative algebra  a couple of algebras which products are respectively commutative and skew-symmetric  and it is linked with  the algebra obtained from the formal deformation. The anti-associative case is developed with a link with the Jacobi-Jordan algebras.



\noindent{\bf Introduction}

In this work $\K$ is a field of characteristic $0$. By a $\K$-algebra, we mean a $\K$-vector space $A$ with a bilinear map $\mu$ and we denote this algebra by $(A,\mu)$. We assume that $\mu$ satisfies a quadratic relation 
%(but this can be extended to a $n$-ary relation) 
denoted $\mu \bullet \mu =0.$ For example, for the associative case, we have $\mu \bullet \mu = \mu \circ ( \mu \otimes Id-Id \otimes \mu)$. The set of $n$-dimensional algebras satisfying a quadratic relation $\mu \bullet \mu =0$ is an algebraic variety $\mathcal{V}_n$ over $\K$ and the classical notion of formal deformation enables a description of  neighborhoods of any point of this variety (for a topology adapted to the structure of the algebraic variety).  A naive definition of a formal deformation of a point $\mu \in \mathcal{V}_n$  is a formal series $\mu_t = \mu + \sum _{k \geq 1} t^k \varphi_k$ considered as a bilinear map on the $\K[[t]]$-algebra $A[[t]]=A \otimes \K[[t]]$,  where the maps $\varphi_k$ are bilinear maps on $A$ which satisfy quadratic relations resulting from the formal identity
$$\mu_t \bullet \mu_t =0.$$
In particular we have $\mu \bullet \mu = 0$ in degree $0$, $\mu \bullet \varphi_1 + \varphi_1 \bullet \mu=0$ in degree $1$, and so on.  Formal deformations are mainly used for the local study of $\mathcal{V}_n$. For example
a point of this variety with only  isomorphic deformations  is topologically rigid, that is, its orbit  is open under the natural action of the linear group. But there are other applications of these deformations. If we consider a formal deformation of a given point $\mu$, it determines new algebra multiplications that are related to the original one. In fact the linear term $\varphi_1$ of $\mu_t$ is also a multiplication on $A$ whose quadratic relation is a consequence of the degree $1$ relation  $\mu \bullet \varphi_1 + \varphi_1 \bullet \mu=0$. A fundamental consequence is  deformation quantization theory  introduced in \cite{Li}. In a simplified way, if we consider a formal associative deformation  $\mu_t=\mu_0+t\varphi_1+\cdots$  of a commutative associative algebra $(A,\mu_0)$, the first term $\varphi_1$ is a cocycle for the Hochschild cohomology associated with $(A,\mu_0)$ and is also a Lie-admissible multiplication whose associated Lie bracket $\psi_1$ satisfies the Leibniz identity with the initial commutative associative multiplication. Then this formal deformation naturally determines  a Poisson algebra $(A,\mu_0,\psi_1)$ and the algebra  $(A[[t]],\mu_t)$ is  a   deformation  quantization of the Poisson algebra  $(A,\mu_0,\psi_1)$. One of the aims of this work is to extend this construction for nonassociative algebras. As the nonassociative algebra world is very wide, we will focus on a class of nonassociative algebras whose quadratic defining relation has symmetric properties linked with the symmetric group and previously studied in \cite{G.R.Nonass}.  

To enlarges the classical deformation  quantization of Poisson algebras,  we have to consider weakly-associative formal deformations of a commutative associative product $\mu_0$, that is formal deformations of $\mu_0$  in the category of weakly-associative algebras as $\varphi_1$ associated with the formal deformation 
$$\mu_t=\mu_0+t\varphi_1+\cdots$$
implies that $(A,\mu_0, \psi_1)$  (where $\psi_1$ is the anticommutator of $\varphi_1$) is still a Poisson algebra. 

We also  consider deformation quantization of generalization of Poisson algebras
and for this consider two generalizations of a  Poisson algebra. Recall that a Poisson algebra is an algebra  $(A,\bullet, \{ \, , \, \})$ with a  commutative associative multiplication $\bullet$ and a Lie bracket  $\{ \, , \, \}$ tied up by the Leibniz identity
$$\mathcal{L}_{\bullet,\{\, , \, \}}(x,y,z)=\{x,y\bullet z\}-y \bullet \{x,z\}-\{x,y\} \bullet  z=0 $$
We already defined in \cite{RWass}   nonassociative Poisson algebras  which are $(A,\bullet, \{ \, , \, \})$  where  $\bullet$ is a commutative multiplication and  $\{ \, , \, \}$ Lie bracket  tied up by the Leibniz rule. So the axioms of a  nonassociative Poisson algebra $\{ \, , \, \}$  are the axioms of a Poisson algebra relaxing the associativity of the operation $\bullet $.

Another generalization is introduced in this paper : the $v$-Poisson algebras where $v$ is a vector of $\KS$. The axioms of a $v$-Poisson algebra are those of  a Poisson algebra weakening the Leibniz identity using the vector $v$. For example  for a $(Id-\tau_{12})$-Poisson algebra the Leibniz rule is replaced by $v$-Leibniz rule
$$\mathcal{L}_{\bullet,\{\, , \, \}}(x,y,z)-\mathcal{L}_{\bullet,\{\, , \, \}}(y,x,z)=0.$$

We can also define a nonassociative $v$-Poisson algebra $(A,\bullet, \{ \, , \, \})$ considering a  commutative multiplication $\bullet$ and a Lie bracket $\{\, , \, \}$ tied up by the $v$-Leibniz rule.

In particular $v$-algebras corresponding to $v=Id+c+c^2$, also called $G_5$-algebras  or $A_3$-associative algebras \cite{Pol}, give the largest class to obtain quantization deformation of a $v$-Poisson algebra by taking $v$-formal deformation of a commutative associative product \cite{RWass}. From a  $G_5$-formal deformations $$\mu_t=\mu_0+t\varphi_1+\cdots$$ of a commutative associative product $\mu_0$ we obtain the algebra $(A,\mu_0,\psi_1)$ which is not a Poisson algebra but a 
$v$-Poisson algebra.

A useful trick to understand the properties of the algebra obtained by formal deformation  or deform a given algebra in a good class  is
%n order to see which algebras are associated with a given commutative algebra, we can also 
to use the  
polarization/depolarization  process introduced in \cite{RM} in the case of Poisson algebras.  Considering  a nonassociative multiplication, this process consists in looking  at the properties of the symmetric and skew-symmetric bilinear applications that are attached to it. We develop in Section 5 the polarization/depolarisation trick for the algebras studied in the first sections. 
A similar study on the  link between polarization and deformations have been done in \cite{Pol}. 

The polarization/depolarization shows that there is a one-to-one correspondence between nonassociative $(Id+c+c^2)$-Poisson algebras and $G_5$-associative algebras (Section 5.5).

The $\KS ^2$-algebra case which is a generalization of the $v$-algebra case including Leibniz algebras is investigated   in Section 4 and we study  $\KS^2$-formal deformations of these algebras. The polarization/depolarization process is also developed in Section 6. A particular look is given to the anti-associative case that is related to the relation $(xy)z+x(yz)=0$. We recall that the corresponding operad is non Kozsul, the description of the "natural cohomology" and the cohomology of the minimal model which parametrizes the deformations \cite{M-R-galgalim}. The deformation quantization process concerns in this case skew-symmetric anti-associative algebras which are related to {\it anti-Poisson algebras} which are defined in Theorem  \ref{antiP} and where the Lie Poisson bracket is replaced by a Jacobi-Jordan product also called a mock-Lie product (see \cite{BF,Camacho}). So we get Jacobi-Jordan algebras by polarization of anti-associative algebras and we study the corresponding operads and describe free Jacobi-Jordan algebras with small number of generators.

This paper also gives  a generalization of the Leibniz identity in a graded version (see Equation (\ref{graded Leibniz})  ) which gives the usual Leibniz identity for $(\rho,\psi)$ a couple of (commutative - skew-symmetric) multiplications, but also Jacobi identity  for $(\psi,\psi)$ with a skew-symmetric multiplication $\psi$. If we consider it for $(\rho,\rho)$ with a commutative multiplication $\rho$, we then obtain the Jacobi-Jordan identity. We also obtain for $(\psi,\rho)$ a couple of (skew-symmetric - commutative) multiplications an identity appearing in the anti-associative algebra case.

\tableofcontents

\section{$\K[\ST]$-associative algebras}
Let  $\Sigma_3=\{Id,\tau_{12},\tau_{13},\tau_{23},c,c^2\}$ be the symmetric group of degree $3$ where $\tau_{ij}$ is the transposition between $i$ and $j$ and $c$ the cycle $(231)$. The product $\sigma\sigma'$  corresponds to the composition $\sigma \circ \sigma'$. Let $\K[\Sigma_3]$ be the group algebra of $\Sigma_3$. It is provided with an associative algebra structure and with a $\Sigma_3$-module structure. The left-action of $\ST$ on  $\KS$ is given by
$$(\sigma \in \ST, v=\sum a_i\sigma_i \in \KS) \ra \sum a_i \sigma  \sigma_i.$$ For any $v \in \KS,$ the corresponding orbit is  denoted by $\mathcal{O}_l(v)=\{v, \tau_{12}v,\! \tau_{13}v, \tau_{23}v,cv,c^2v\}$ or simply $\mathcal{O}(v)$ and $F_v={\rm Span}(\mathcal{O}(v))$  is the $\K$-linear subspace of $\KS$ generated by $\mathcal{O}(v)$. It is also a $\ST$-module. 

\noindent Some notations:  (1) We call  canonical basis of $\KS$ the ordered family $\! \{Id,\tau_{12},\tau_{13},\tau_{23},c,c^2\! \}$ and $(a_1,a_2,a_3,a_4,a_5,a_6)$ are the coordinates of the vector $v=a_1Id+a_2\tau_{12}+a_3\tau_{13}+a_4\tau_{23}+a_5 c+a_6 c^2$ in the canonical basis. We denote $M_v$ the matrix composed of the column component vectors of the family $( v,\! \tau_{12}v,\! \tau_{13}v,\! \tau_{23}v,cv,\! c^2v)$ in the canonical basis:
$$
M_v=
\left(
\begin{array}{cccccc}
a_1 & a_2 & a_3 & a_4 & a_6 & a_5 \\
a_2	 & a_1 & a_5 & a_6 & a_4 & a_3 \\
a_3 & a_6 & a_1 & a_5 & a_2 & a_4 \\
a_4 & a_5 & a_6 & a_1 & a_3 & a_2 \\
a_5 & a_4 & a_2 & a_3 & a_1 & a_6 \\
a_6 & a_3 & a_4 & a_2 & a_5 & a_1 \\
\end{array}
\right)
$$

\noindent (2) Let $A$ be a $\K$-vector space. The symmetric monoidal structure on the category of vector spaces naturally turns $A^{\otimes 3}$ in a representation of $\Sigma_3$. We denote this representation by
$$\Phi^A: \Sigma_3 \rightarrow Aut(A^{\otimes 3}) \subset End(A^{\otimes 3}).$$
The universal property of the group algebra allows one to extend this representation to
$$\Phi^A: \KS \rightarrow End(A^{\otimes 3}).$$
Thus, if  $X,Y,Z$ are three vectors of $A$, and if we denote $\Phi^A_\sigma $ instead of $\Phi^A(\sigma)$, we have:
$$\begin{array}{   l | r | r |  r | r | r | r}
 & (X,Y,Z) &(Y,X,Z) & (Z,Y,X) &  (X,Z,Y) & (Y,Z,X) & (Z,X,Y) \\
 \hline
\Phi^A_{Id}            & (X,Y,Z) &(Y,X,Z) & (Z,Y,X) &  (X,Z,Y) & (Y,Z,X) & (Z,X,Y) \\
\hline
\Phi^A_{\tau_{12}} & (Y,X,Z) & (X,Y,Z) & (Y,Z,X) & (Z,X,Y) & (Z,Y,X) & (X,Z,Y) \\
\hline 
\Phi^A_{\tau_{13}} & (Z,Y,X) & (Z,X,Y) & (X,Y,Z) & (Y,Z,X) & (X,Z,Y) & (Y,X,Z) \\
\hline 
\Phi^A_{\tau_{23}} & (X,Z,Y) & (Y,Z,X) & (Z,X,Y) & (X,Y,Z) & (Y,X,Z) & (Z,Y,X) \\
\hline 
\Phi^A_{c} & (Y,Z,X) &(X,Z,Y) & (Y,X,Z) &  (Z,Y,X) & (Z,X,Y) & (X,Y,Z) \\
\hline 
\Phi^A_{c^2} & (Z,X,Y) &(Z,Y,X) & (X,Z,Y) &  (Y,X,Z) & (X,Y,Z) & (Y,Z,X) \\
\end{array}$$
For any $w =\sum_{i=1}^6 w_i \sigma_i \in \KS$, 
$$\Phi^A_w=\sum_{i=1}^6 w_i \Phi^A_{\sigma_i }.$$
In particular, for any $v,w \in \KS$, 
$$\Phi^A_w \circ \Phi^A_v=\Phi^A_{v  w}.$$

\begin{definition}
Consider a nonzero vector $v \in \mathbb{K}\, [\Sigma_3]$.
An algebra $(A,\mu)$ is a $v$-associative algebra or simply a $v$-algebra if 
$$ \mathcal{A}_{\mu} \circ \Phi^A_v =0,
$$
 where $\mathcal{A}_\mu$ is the associator of $\mu$ that is $\mathcal{A}_\mu(X,Y,Z)=\mu(\mu(X,Y),Z)-\mu(X,\mu(Y,Z)).$
 
We will also say that the algebra $A$ is $v$-associative.
\end{definition}
A  $v$-algebra $(A,\mu)$ is  also a $v'$-algebra if $v' \in F_{v}={\rm Span}(\mathcal{O}(v))$. But for any $v' \in F_v$ such that $\dim F_{v'} < \dim F_v$, the $v'$-associativity doesn't imply the $v$-associativity. For example, if $v=Id-\tau_{12} +c$, the vector $v'=Id+\tau_{13} $ is in $F_v$. But $\dim F_{v'}=3 < \dim F_v=4$ and the $(Id+\tau_{13})$-associativity doesn't imply the $(Id-\tau_{12} +c)$-associativity. Of course $v$-algebras are the same that $\sigma(v)$-algebras for any $\sigma \in \Sigma_3$ as we trivially have that 
$F_v=F_{\sigma(v)}$ and more generally $v$-algebras are the same as $v'$-algebras if and only if $F_v=F_{v'}$

\medskip

We obtain from the $\ST$-module structure of $F_v$  that it decomposes in a direct sum of $F_{v_i}$ associated with the irreducible representations of $\ST$.
There exist two particular vectors in $\KS$ denoted here by $v_{Lad}$ and $v_{3Pa}$ 
%(we will explain the notation meaning later) 
corresponding to the only $1$-dimensional irreducible signum  and  trivial representations:
$$v_{Lad}=\ds \sum_{\sigma \in \ST} \varepsilon(\sigma)\sigma \quad {\rm and} \quad v_{3Pa}=\ds \sum_{\sigma \in \ST} \sigma,$$
where $\varepsilon(\sigma)$ is the signature of the permutation $\sigma$. The vectors $v_{Lad}$ and $v_{3Pa}$ are the unique vectors $v$ such that $F_v$ is one dimensional up to a scalar factor. 

\begin{proposition}
An algebra $(A,\mu)$ is 
\begin{enumerate} 
\item Lie-admissible if and only if it is $v_{Lad}$-associative,
\item  $3$-power-associative if and only if it is $v_{3Pa}$-associative.
\end{enumerate} 
\end{proposition}
\pf See \cite{G.R.Nonass}. 
The classes of Lie-admissible algebras and power-associative algebras have been introduced by Albert in \cite{Al}.  An algebra is called Lie-admissible if the skew-symmetric  bilinear map $\psi$ related to $\mu$ is a Lie bracket. This is equivalent to write $\mathcal{A}_\mu \circ \Phi^A_{v_{Lad}}=0$. An algebra  is said to be power-associative if the subalgebra generated by any element is associative. 
Over a field of characteristic $0$, an algebra is power-associative if  it satisfies $\mathcal{A}_{\mu}(x,x,x)=\mathcal{A}_{\mu}(x^2,x,x)=0$
for any $x \in A$. An algebra  is said to be $3$-power-associative  if and only if it satisfies $\mathcal{A}_{\mu}(x,x,x)=0$ for any $x \in A$. This last condition is equivalent, by linearization, to $\mathcal{A}_\mu \circ \Phi^A_{v_{3Pa}}=0.$

\medskip

\noindent{\bf Remark.} If $(A,\mu)$ is  $3$-power-associative, then  $\mathcal{A}_\mu \circ \Phi^A_{v_{3Pa}}=0$ which implies
$\mathcal{A}_\mu(x^2,x,x)+\mathcal{A}_\mu(x,x^2,x)+\mathcal{A}_\mu(x,x,x^2)=0.$ But we also have $\mathcal{A}_\mu(x^2,x,x)-\mathcal{A}_\mu(x,x^2,x)+\mathcal{A}_\mu(x,x,x^2)=0.$ In fact, since $xx^2=x^2x$, then 
$$\mathcal{A}_\mu(x^2,x,x)-\mathcal{A}_\mu(x,x^2,x)+\mathcal{A}_\mu(x,x,x^2)= x^2(x^2)-(x^2x)x-x(x^2x)+(xx^2)x+x(xx^2)-x^2x^2=0.$$
We deduce $\mathcal{A}_\mu(x,x^2,x)=0$ and  $\mathcal{A}_\mu(x^2,x,x)+\mathcal{A}_\mu(x,x,x^2)=0.$ Then  a $3$-power-associative algebra is  power-associative if and only if  
$\mathcal{A}_\mu(x^2,x,x)-\mathcal{A}_\mu(x,x,x^2)=0.$ 
So a sufficient condition for a $3$-power-associative to be power-associative is $\mathcal{A}_\mu \circ \Phi_{Id-\tau_{13}}= 0.$

\medskip

There is a third irreducible representation of the group $\ST$, the first two being associated with the vectors $v_{Lad}$ and $v_{3Pa}$. It is a representation of degree $2$. It will be used   later when we give the classification of $v$-algebras which are Lie-admissible or $3$-power associative algebras using  the rank of $v$.

\section{Formal deformations of $v$-algebras}
\subsection{Generalities}
Let $(A,\mu_0)$ be  a $v$-algebra where  $v=\sum a_i\sigma_i \in \KS$.  Let $\varphi_1$ and $\varphi_2$ be two bilinear maps on $A$. We denote by 
$\varphi_1 \bullet \varphi_2$ the trilinear map on $A$
$$\varphi_1 \bullet \varphi_2=\varphi_1 \circ (\varphi_2 \otimes Id)-\varphi_1 \circ (Id \otimes \varphi_2)$$
and
$\varphi_1 \bullet_v \varphi_2$ the trilinear map on $A$
%$$\varphi_1 \bullet_v \varphi_2 (x_1,x_2,x_3)=\sum a_i\left( \varphi_1(x_{\sigma_i(1)},\varphi_2(x_{\sigma_i(2)},x_{\sigma_i(3)}))-\varphi_1(\varphi_2(x_{\sigma_i(1)},x_{\sigma_i(2)}),x_{\sigma_i(3)})\right)$$
%for any $x_1,x_2,x_3 \in A,$
%that is,
$$\varphi_1 \bullet_v \varphi_2=(\varphi_1 \circ (\varphi_2 \otimes Id)-\varphi_1 \circ (Id \otimes \varphi_2)) \circ \Phi_v=(\varphi_1 \bullet \varphi_2) \circ \Phi_v.$$

Let $(A,\mu_0)$ be a $v$-algebra. A $v$-formal deformation of $(A,\mu_0)$ is given by a family of bilinear maps on $A$
$$\{\varphi_j : A \otimes A \rightarrow A , \  j\in \N\}$$
with $\varphi_0=\mu_0$ and
satisfying
\begin{equation}
\label{def}
\sum_{i+j=k, i,j \geq 0}\varphi_i \bullet_v \varphi_j=0, \ \ k\geq 0.
\end{equation}
If we denote by $\K[[t]]$ the algebra of formal series with one indeterminate $t$, this definition is equivalent  to consider on the space $A [[t]] =A \otimes \K[[t]]$ (recall that $A$ is of finite dimensional) of formal series with coefficients in $A$ a structure  of $\K[[t]]$-$v$-associative algebra such that the canonical map $A[[t]]/ tA[[t]] \rightarrow A$ is an isomorphism of $v$-algebras. It is useful to write
$$\mu_t=\mu_0+t \varphi_1+t^2\varphi_2+ \cdots$$
Equation (\ref{def}) implies at the order $k=0$ that $\mu_0$ is $v$-associative. The order $k=1$ writes
$$\mu_0 \bullet_v \varphi_1 + \varphi_1 \bullet_v \mu_0=0.$$
To be consistent with the conventional cohomological approaches to deformations, we will denote by 
$\delta_{v,\mu_0}^2 \varphi$ the trilinear map
$$\delta_{v,\mu_0}^2 \varphi=\mu_0 \bullet_v \varphi + \varphi \bullet_v \mu_0.$$
In fact, we know  that a cohomological complex which parametrizes formal deformations of  algebras  over a quadratic operad exists and $\delta_{v,\mu_0}^2$ corresponds to the second coboundary operator. For example, if $v=Id$, then $(A,\mu_0)$ is associative and $\delta_{Id,\mu_0}^2$ is the coboundary operator associated with the Hochschild complex of $A$ classically denoted by
$\delta_{H,\mu_0}^2$ and we have
$$\delta_{H,\mu_0}^2\varphi(x,y,z)=-x\varphi(y,z)+\varphi(xy,z)-\varphi(x,yz)+\varphi(x,y)z=(\mu_0\bullet \varphi + \varphi \bullet \mu_0)(x,y,z)$$
where, to simplify the notations, $xy$ means $\mu_0(x,y)$. Then
$$\delta_{v,\mu_0}^2 \varphi=\delta_{H,\mu_0}^2 \varphi \circ \Phi_v.$$Coming back to Equation (\ref{def}), we obtain
\begin{equation}\label{eq}
\begin{array}{ll}
\medskip
 {\rm order} \  0 :& \mathcal{A}_{\mu_0} \circ \Phi_v =0,      \\
 \medskip
   {\rm order} \  1 :  &   \delta_{v,\mu_0 }^2\varphi_1 =0,\\
   \medskip
   {\rm order} \  2 : & \varphi_1 \bullet_v \varphi_1 +  \delta_{v,\mu_0}^2 \varphi_2 =0.
\end{array}
\end{equation}
Let  $v_1$ be a vector in $\KS$. Then $\varphi_1 \bullet_v \varphi_1 \circ \Phi_{v_1}=0$, that is $(A, \varphi_1)$ is a $v_1v$-algebra,  if and only if 
$\delta_{v,\mu_0}^2 \varphi_2 \circ \Phi_{v_1}=0.$ But $\delta_{v,\mu_0}^2 \varphi_2 \circ \Phi_{v_1}=\delta_{H,\mu_0}^2  \varphi_2 \circ \Phi_{v_1v}.$  So we will look when $\delta_{H,\mu_0}^2  \varphi_2 \circ \Phi_{v_1v}=0$ is satisfied but asking moreover the commutativity of the multiplication $\mu_0$. This new hypothesis will be justified in the study of  deformation quantization in the next section.

%%%%%%%%%%%%%%%
\subsection{Case of a commutative $v$-algebra $(A,\mu_0)$}
\begin{lemma}\label{1}
Let $(A,\mu_0)$ be a commutative algebra with $\mu_0 \neq 0$ and $\delta^2_{H,\mu_0}$  the Hochschild coboundary operator:
$$\delta^2_{H,\mu_0} \varphi (X,Y,Z)=-X\varphi(Y,Z)+\varphi(XY,Z)-\varphi(X,YZ)+\varphi(X,Y)Z$$
where $X,Y,Z\in A$, the map $\varphi$ is bilinear on $A$  and  $XY$denotes the product $\mu_0(X,Y)$.  
Then $ \delta^2_{H,\mu_0}\varphi \circ \Phi_{v_{Lad}}= 0.$
\end{lemma}
\pf It is easy to see that 
$$\begin{array}{rl}
\delta^2_{H,\mu_0}  \varphi  \circ \Phi_{v_{Lad}} (X_1,X_2,X_3)=& \ds \sum_{\sigma \in \mathbb{K}[\Sigma_3]} \varepsilon(\sigma)   \delta^2_{H,\mu_0}\varphi (X_{\sigma(1)},X_{\sigma(2)},X_{\sigma(3)})\\
=&\ds \sum_{\sigma \in \mathbb{K}[\Sigma_3]} \varepsilon(\sigma)    (\varphi(X_{\sigma(1)},X_{\sigma(2)})X_{\sigma(3)}-X_{\sigma(1)}\varphi(X_{\sigma(2)},X_{\sigma(3)}))
\\ & \ds +\sum_{\sigma \in \mathbb{K}[\Sigma_3]} \varepsilon(\sigma)   \varphi(X_{\sigma(1)}X_{\sigma(2)},X_{\sigma(3)}) \\
& \ds -\sum_{\sigma \in \mathbb{K}[\Sigma_3]} \varepsilon(\sigma) \varphi(X_{\sigma(1)},X_{\sigma(2)}X_{\sigma(3)})
\end{array}$$ 

 The commutativity  of $\mu_0$ implies the cancellation of each term.  $\Box$

\smallskip

\noindent Let us apply Lemma \ref{1} to study, for a commutative $v$-associative  algebra $(A,\mu_0)$, the equation
$$\varphi_1 \bullet_v \varphi_1 +  \delta_{v,\mu_0}^2 \varphi_2 =0$$
For any $v_1 \in \KS$ we have
$$\varphi_1 \bullet_v \varphi_1\circ \Phi_{v_1} +  \delta_{v,\mu_0}^2 \varphi_2\Phi_{v_1}=\varphi_1 \bullet_{v_1v} \varphi_1 +  \delta_{v_1v,\mu_0}^2 \varphi_2 =0.$$
There is an obvious solution to this equation corresponding to the case $v_1 v=0$ but which do not lead to any properties on $\varphi_1.$  If $M_v$ is the matrix associated with $v$, then the equation $v_1v=0$ corresponds to the linear system $M_v V_1=0$ where $V_1$ is the column matrix of the vector $v_1$ and $v_1v=0$ if and only if $V_1 \in \ker  M_v$. 
The rank of $M_v$ is the dimension of $Span(\mathcal{O}(v))=F_v$ and it is maximal if and only if $Id\in F_v$. In this case $\mu_0$ is associative. In all the other cases, $\rk(M_v) < 6$ and $\dim \ker M_v \geq 1$.  For example, if $v=v_{Lad}$, then $\dim F_v=1$ and $v_1v_{Lad} \in F_{v_{Lad}}$ for any $v_1 \in \KS$. More precisely if $v_1=(a_1,a_2,a_3,a_4,a_5,a_6)\in \KS$ then $v_1v_{Lad}=\lambda v_{Lad} $ with $\lambda=a_1-a_2-a_3-a_4+a_5+a_6$ and 
$\{ v_1 \in \KS; v_1v_{Lad}=0\}$ is a $5$-dimensional subspace of $\KS$ defined by the linear equation $a_1-a_2-a_3-a_4+a_5+a_6=0$.

\noindent  More generally, let $v\in \KS$. We consider the equation $v_1v=v_{Lad}$, which corresponds to the linear system 
$M_v V_1=V_{Lad}$. If $v_1 \notin \ker M_v$, the equation $v_1v=v_{Lad}$ is equivalent to $v_{Lad} \in \im M_v.$.

\begin{proposition}\label{Mv}
Let $(A,\mu_0)$ be a commutative $v$-algebra and $\mu_t=\mu_0+t\varphi_1+t^2\varphi_2+ \cdots $ be a $v$-formal deformation of $\mu_0$. Then $\varphi_1$ is a Lie-admissible multiplication on $A$ if  $v_{Lad} \in \im M_v$ that is $v_{Lad} \in F_v$.
\end{proposition}
\pf If $v_{Lad} \in F_v$, then there exists $v_1$ such that $v_1v=v_{Lad}$ and 
$$\varphi_1 \bullet_{v_1v} \varphi_1 +  \delta_{v_1v,\mu_0}^2 \varphi_2 =\varphi_1 \bullet_{v_{Lad}} \varphi_1 +  \delta_{v_{Lad},\mu_0}^2 \varphi_2 =0.$$
From Lemma \ref{1}, we deduce
$$\varphi_1 \bullet_{v_{Lad}} \varphi_1 =0.$$

\begin{theorem}\label{co}
Let $w=a_1Id+a_2\tau_{12}+a_3\tau_{13}+a_4\tau_{23}+a_5c+a_6c^2$ be a vector of $\KS$. 
\begin{enumerate}
\item $a_1-a_2-a_3-a_4+a_5+a_6 \neq 0$ if and only if 
$v_{Lad} \in F_w$.
\item 
$a_1+a_2+a_3+a_4+a_5+a_6 \neq 0$ if and only if 
$v_{3Pa} \in  F_w$.
\end{enumerate}
\end{theorem}
\pf In fact,  if $\lambda=a_1-a_2-a_3-a_4+a_5+a_6 \neq 0$,  then $\lambda$ is a non null eigenvalue of $M_w$ and $v_{Lad}$ is an eigenvector corresponding to $\lambda$. In this case $v_{Lad} \in \im M_w$ that is  $v_{Lad} \in  F_w$.
Suppose now that $v_{Lad}\in F_w$. Then it implies that there exists  $ v=\lambda_1 Id+\lambda_2 \tau_{12}+ \lambda_3  \tau_{13}+\lambda_4  \tau_{23}+\lambda_5 c +\lambda_6 c^2$ such that   $v w= v_{Lad}$ that is  
$$\begin{array}{rl}
 M_w V=\ ^t(1,-1,-1,-1,1,1)
 \Leftrightarrow 
 \begin{array}{l}
 L_1\\
 L_2\\
 L_3\\
 L_4\\
 L_5\\
 L_6
 \end{array}
 \!
   \left\{
 \begin{array}{l} 
 a_1\lambda_1+a_2\lambda_2+a_3\lambda_3+a_4\lambda_4+a_6 \lambda_5+a_5\lambda_6=1\\
 a_2\lambda_1+a_1\lambda_2+a_5\lambda_3+a_6\lambda_4+a_4 \lambda_5+a_3\lambda_6=-1\\
 a_3\lambda_1+a_6\lambda_2+a_1\lambda_3+a_5\lambda_4+a_2 \lambda_5+a_4\lambda_6=-1\\
 a_4\lambda_1+a_5\lambda_2+a_6\lambda_3+a_1\lambda_4+a_3 \lambda_5+a_2\lambda_6=-1\\
 a_5\lambda_1+a_4\lambda_2+a_2\lambda_3+a_3\lambda_4+a_1 \lambda_5+a_6\lambda_6=1\\
a_6\lambda_1+a_3\lambda_2+a_4\lambda_3+a_2\lambda_4+a_5 \lambda_5+a_1\lambda_6=1\\
 \end{array}
 \right. .
  \end{array}
$$
Consider  $L_1-L_2-L_3-L_4+L_5+L_6$ the system  implies that $(a_1-a_2-a_3-a_4+a_5+a_6)(\lambda_1 -\lambda_2 - \lambda_3 - \lambda_4 + \lambda_5  +\lambda_6)=6$, that is, $a_1-a_2-a_3-a_4+a_5+a_6\neq 0$.

We apply the same technique for the 3-power associative case to prove that  $a_1+a_2+a_3+a_4+a_5+a_6 \neq 0$ if and only if 
$v_{3Pa} \in  F_w.$ $\square$

%For the general case, in \cite{G.R.Nonass}, we determine vectors $v$ such that $v_{Lad} \in F_v$. 
\medskip Note that this result can be interpreted in terms of representations of $\KS$. The projections to the isotopic components of the regular reprecentation of $\KS$ are given by right multiplication by the Young symmetrizers. In the case of the trivial and sign component, the Young symmetrizers are proportional to $v_{3Pa}$ and $v_{Lad}$ so that the subrepresentation $F_w$ generated by a vector $w$ contains a copy of the sign representation $F_{v_{Lad}}$ if and only if $wv_{Lad} \neq 0.$

\medskip

\noindent{\bf Remarks.} Consider a commutative product $\mu_0$ 
% we want to know what are the conditions on  $w$ such that   $ \delta^2_{H,\mu_0}\varphi \circ \Phi_w= 0$ for any bilinear map $\varphi$. We already now from Lemma \ref{1} that $w=v_{Lad}$ is fine.
%Consider 
and $w=a_1Id+a_2\tau_{12}
+a_3\tau_{13}+a_4\tau_{23}+a_5c+a_6c^2$  a vector of $\KS$. The commutativity of  $\mu_0$ implies \begin{equation}
\label{ww}
\begin{array}{l}
\delta^2_{H,\mu_0}\varphi \circ \Phi_w(X,Y,Z)= \\

(a_5-a_1)X\varphi(Y,Z) +(a_3-a_4)X\varphi(Z,Y)+(a_4-a_2)Y\varphi(X,Z)+(a_6-a_5)Y\varphi(Z,X)\\
+(a_1-a_6)Z\varphi(X,Y)+(a_2-a_3)Z\varphi(Y,X)+(a_1+a_2)\varphi(XY,Z)+(a_3+a_5)\varphi(YZ,X)\\
+(a_4+a_6)\varphi(XZ,Y)-(a_1+a_4)\varphi(X,YZ)-(a_2+a_5)\varphi(Y,XZ)-(a_3+a_6)\varphi(Z,XY).
\end{array}
\end{equation}
Remark that   we reobtain  Lemma \ref{1} considering $a_1=-a_2=-a_3=-a_4=a_5=a_6$.

If $\dim A=1$, that is $A=\K$, since $\mu_0 \neq 0$, any bilinear form writes $\varphi=a\mu_0$ with $a \in \K$ and satisfies 
$\delta^2_{H,\mu_0}\varphi \circ \Phi_w=0$ with no assumption on $w.$

Suppose that $\dim A \geq 2$.  Equation (\ref{ww}) with $X=Y=Z$ reduces to
$$(a_1+a_2+a_3+a_4+a_5+a_6)(\varphi(XX,X)-\varphi(X,XX))=0.$$
Since $A^2 \neq 0$, we can choose $X$  such that $XX\neq 0$.  If $\varphi$ is such  that  $\varphi(XX,X)-\varphi(X,XX)\neq 0$ we obtain
$$a_1+a_2+a_3+a_4+a_5+a_6=0$$
and $w=a_1Id+a_2\tau_{12}
+a_3\tau_{13}+a_4\tau_{23}+a_5c-(a_1+a_2+a_3+a_4+a_5)c^2$.  
Then from Theorem \ref{co}, $v_{3Pa} \notin F_{w}$.  For example if $\varphi$ is skew-symmetric with  $\varphi(XX,X)\neq 0$ for some $X \in A$, the fact that $\delta^2_{H,\mu_0}\varphi \circ \Phi_w=0$  implies that $v_{3Pa} \notin F_{w}$.

\section{Deformation quantization of the $v$-algebras with $v_{Lad} \in F_v$}

Recall that  the rank of a vector $v \in \KS$ is the dimension of the vector space $F_v={\rm Span}(\mathcal{O}(v))$. Then, if $v \neq 0$, then $1 \leq \rk(v) \leq 6.$ If $\rk(v)=6$ we have $F_v=\KS$ and $Id \in F_v$. In this case, any $v$-associative algebra is associative and we can assume that $v=Id$. Similarly, if $\rk(v)=1$, then $\dim F_v=1$ and it is  a one-dimensional invariant subspace of $\KS$. We have seen that, in this case, $v=v_{Lad}=Id-\tau_{12}-\tau_{13}-\tau_{23}+c+c^2$ or $v=v_{3Pa}=Id+\tau_{12}+\tau_{13}+\tau_{23}+c+c^2$. In this section we will focus on $v$-algebras such as $v_{Lad} \in F_v$ because of Proposition \ref{Mv}.  In  \cite{G.R.Nonass}, we have the following result:
\begin{theorem}
Every Lie-admissible $v$-algebra $(A,\mu_0)$ corresponds to one of the following types:
\begin{enumerate}
\item Type (I): $dim \ F_v=1$ and $F_v=F_{v_{Lad}}$.

\item Type (II): $dim \ F_v=2$ and $v= Id+c+c^2$. The corresponding $v$-algebras are also  $3$-power-associative algebras. These algebras correspond to the $G_5$-associative algebras.

\item Type (III): $dim \ F_v=3$. The corresponding  $v$-algebras satisfy:
$$v=\alpha Id-\alpha \tau_{12} +(\alpha+\beta-3)\tau_{13}-\beta \tau_{23}+\beta c +(3-\alpha -\beta) c^2$$
with $(\alpha,\beta) \neq (1,1)$.
\item Type (IV): $dim \ F_v=4.$ The $v$-algebras are of the following type:
\begin{enumerate}
\item  (IV.1): $v=2Id +(1+t)\tau_{12} +\tau_{13} +c+ (1-t) c^2$ with $t \neq 1$,
\item (IV.2): $v=2Id +\tau_{12} +\tau_{23} +c+c^2$ 
\end{enumerate}
\item Type (V): $dim \ F_v=5$ and $v= 2Id -\tau_{12}-\tau_{13}-\tau_{23}+c$
\item Type (VI): $dim \ F_v=6.$ This corresponds to the class of associative algebras that is $v=Id$.
\end{enumerate}
\end{theorem}

The $v$-algebras associated with the vector $v=Id+c+c^2$ have  $\mathcal{O}(v)=\{Id+c+c^2, \tau_{12}+\tau_{13}+\tau_{23}\}$ and $\dim F_v=2.$ For example, $(A,\mu)$ with $\mu$ skew-symmetric is a $v$-algebra if and only if it is a Lie algebra. 

In \cite{G.R.LieAdm} we have studied 
particular classes of $v$-algebras called 
$G$-associative algebras whose defining quadratic relation  is associated with  subgroups of $\Sigma_3$. Consider $G_1=\{Id\},G_2=\{Id,\tau_{12}\},G_3=\{Id,\tau_{13}\},G_4=\{Id,\tau_{23}\},G_5=A_3=\{Id,c,c^2\}$ and $G_6=\Sigma_3$ the subgroups of $\Sigma_3$. A $G_i$-associative algebra is defined by the relation given by $v_i$-associative algebra with
$$v_i=\sum_{\sigma \in G_i} \varepsilon(\sigma)\sigma.$$
In particular $G_1$-associative algebras correspond to the associative algebras, $G_6$-associative algebras to 
the Lie-admissible algebras. These algebras, as well as the case $G_5=A_3$, have been studied previously. The remaining cases are associated with a vector $v$ of rank 3:
the $G_2$-associative algebras also called Vinberg algebras and associated with the vector $v=Id-\tau_{12}$ correspond to $\alpha=3,\beta=0$; the $G_3$-associative algebras also called Pre-Lie algebras and associated with the vector  $v=Id-\tau_{13}$ correspond to $\alpha=0,\beta=0$ and finally, 
the $G_4$-associative algebras, associated with the vector  $v=Id-\tau_{23}$, correspond to $\alpha=0,\beta=3$.  
We begin this study by the more classical case corresponding to an associative and commutative multiplication $\mu_0$. 
 
 %%%%%%%%%%%%%%%%%%%%%%
\subsection{Rank$(v)=6$: $v=Id$. The associative case}
The study of deformations of associative algebras was initiated by Gerstenhaber \cite{Ge} and deformation quantization by Bayen, Flato, Fronsdal, Lichnerowicz and Sternheimer in \cite{Li}. In a first step,  we summarize this study as part of the $v$-associative algebras.

When $v=Id$, a $v$-algebra is an associative algebra. Let $\mu=\mu_0+t\varphi_1+t^2\varphi_2 + \cdots$ be an associative  formal deformation of a commutative associative multiplication $\mu_0$.  In this case $\delta_{Id,\mu_0}^2=\delta_{H, \mu_0}^2$ and  Equations (\ref{eq}) write
$$\left\{\begin{array}{ll}
 \medskip
   {\rm order} \  1   &   \delta_{H, \mu_0 }^2\varphi_1 =0,\\
   \medskip
   {\rm order} \  2  & \varphi_1 \bullet \varphi_1 +  \delta_{H, \mu_0}^2 \varphi_2 =0.
\end{array}
\right.
$$
From Proposition $\ref{Mv},$ as $v_{Lad} \in F_{Id}$, $\varphi_1$ is a Lie-admissible  multiplication. 
%We have seen that  $\delta_{H,\mu_0}^2 \varphi_2 \circ \Phi_{v_{Lad}}=0$ and the previous equation implies 
%$$\varphi_1 \bullet \varphi_1 \circ \Phi_{v_{Lad}} =\mathcal{A}_{\varphi_1}\circ \Phi_{v_{Lad}}=0.$$
%We deduce that $\varphi_1$ is Lie-admissible. 

%\begin{proposition} Let $\mu=\mu_0+t\varphi_1+t^2\varphi_2 + \cdots$ be an associative  formal deformation of a commutative associative multiplication $\mu_0$.Then $\varphi_1$ is a Lie-admissible  multiplication. 
%\end{proposition}

The bilinear map $\varphi_1$ also satisfies $\delta_{H,\mu_0 }^2\varphi_1 =0$ and so $\delta_{H,\mu_0 }^2\varphi_1\circ \Phi_v =0$ for any $v \in \KS$. Let us determine a vector such that this relation involves a relation on the skew-bilinear map $\psi_1$ attached to $\varphi_1$, that is $\psi_1(x,y)=\varphi_1(x,y)-\varphi_1(y,x)$. If $v=a_1Id+a_2\tau_{12}+a_3\tau_{13}+a_4\tau_{23}+a_5c+a_6c^2$ 
with
$$a_5=a_1-a_3+a_4, \ \ a_6=a_1+a_2-a_3,$$
then, writing $xy$ for $\mu_0(x,y)$, we have  $xy=yx$ and $\delta_{H }^2\varphi_1\circ \Phi_v =0$ is equivalent to
$$
\begin{array}{l}
  a_1(\psi_1(xy,z)+\psi_1(xz,y)+\psi_1(zy,x))+a_2(y\psi_1(x,z)+z\psi_1(x,y)-\psi_1(xz,y) \\
 -\psi_1(xy,z))
 +a_3(\psi_1(xz,y)+x\psi_1(y,z)-z\psi_1(x,y))
 +a_4(-x\psi_1(y,z) -y\psi_1(x,z) 
 -\psi_1(xz,y)\\-\psi_1(zy,x))=0\\
\end{array}
$$
for any $a_1,a_2,a_3,a_4 \in \K$. This is equivalent to
$$
\left\{
\begin{array}{l}
\psi_1(xy,z)+\psi_1(xz,y)+\psi_1(zy,x)=0,\\
y\psi_1(x,z)+z\psi_1(x,y)-\psi_1(xz,y) -\psi_1(xy,z)=0,\\
\psi_1(xz,y)+x\psi_1(y,z)-z\psi_1(x,y)=0,\\
-x\psi_1(y,z) -y\psi_1(x,z) 
 -\psi_1(xz,y)-\psi_1(zy,x)=0.\\
 \end{array}
 \right.
 $$
The third identity is the Leibniz identity between the Lie bracket $\psi_1$ and the commutative associative multiplication $\mu_0$. Since the other identities are consequence of the Leibniz identity, we find the classical result
\begin{proposition} If $\mu_t=\mu_0+t\varphi_1+t^2\varphi_2 +\cdots$ is an associative formal deformation of the commutative associative multiplication  $\mu_0$ on $A$, then $(A,\mu_0,\psi_1)$ is a Poisson algebra and the formal deformation $(A[[t]],\mu_t)$ is a deformation quantization of this Poisson algebra.
\end{proposition}
In this proposition, we see that any associative deformation of the commutative associative algebra $(A,\mu_0)$ gives a quantization. But, are there $v$-formal deformations of $\mu_0$ with $v \in \KS$ but $v \notin \ST$ which define a deformation quantization of a Poisson algebra $(A,\mu_0,\psi)$ for some Lie bracket $\psi$?  In \cite{RWass}, we show that there exists a class of nonassociative algebras, called weakly associative algebras, corresponding to the vector $v=Id-\tau_{12}+c$ that answers the previous question: any $v$-formal deformation of a commutative associative algebra defines a deformation quantization of a Poisson algebra. Since the vector $v=Id-\tau_{12}+c$, associated with the weakly associative algebra 
is of rank $4$,  we will briefly rescall this study in the paragraph dedicated to rank $4$.

\subsection{Rank$(v)=1$: $F_v=F_{v_{Lad}}$, that is the  Lie-admissible  algebras}
Let $(A,\mu_0)$ be a commutative $v_{Lad}$-algebra. Let  $\mu_t=\mu_0+ \sum t^i\varphi_i$ be a $v_{Lad}$-deformation of $\mu_0$. Remark that a commutative product is always Lie-admissible. From Lemma \ref{1}, since $\mu_0$ is commutative, for any bilinear map $\varphi$ we have $\delta^2_{H,\mu_0}\varphi \circ \Phi_{v_{Lad}}=0$. This implies that, for any $i \geq 1$, $\delta^2_{v_{Lad},\mu_0} \varphi_i =0.$ In particular $\varphi_1$ is Lie-admissible. 

\begin{proposition}
Let $(A,\mu_0)$ be a commutative Lie-admissible algebra. For any bilinear map $\varphi$ we have
$$\delta^2_{v_{Lad},\mu_0} \varphi=0.$$
If $\mu_t=\mu_0+t\varphi_1+t^2\varphi_2 +\cdots$ is a Lie-admissible formal deformation of $\mu_0$ then the algebra $(A,\varphi_1)$  is a Lie-admissible.
\end{proposition}

Let us note that if $\mu_0$ is a skew-symmetric Lie-admissible multiplication (non necessarily commutative), that is  $\mu_0$ is a Lie bracket, then 
$$\delta^2_{v_{Lad},\mu_0}\varphi=2\delta^2_{CE,\mu_0}\psi$$
where $\psi$ is the skew-symmetric bilinear map associated to $\varphi$ and $\delta^2_{CE,\mu_0}$ the coboundary operator of the Chevalley Eilenberg cohomology of the Lie algebra $(A,\mu_0)$.

\subsection{Rank$(v)=2:  F_v= F_{Id+c+c^2}$ that is $G_5$-algebra  or $A_3$-associative algebra.}
Recall that for any $w \in \K[\Sigma_3]$,  $F_w$ is a $\Sigma_3$-invariant vector space so a direct sum of irreducible vector spaces. The irreducible vector spaces are $1$-dimensional, that is $F_{v_{Lad}}$ and $F_{v_{3Pa}}$ or $2$-dimensional, that is $F_{\alpha(Id-\tau_{12})+\beta(c-\tau_{23})+(\alpha+\beta)(\tau_{13}-c^2)}$ with $\alpha,\beta \in \mathbb{R}.$ As we considered that $v_{Lad} \in F_v$ we have  that $F_v=F_{v_{Lad}} \oplus F_{v_{3Pa}}$ and we can assume that $v= Id+c+c^2$.

Any commutative multiplication $\mu_0$ satisfies $\mu_0 \bullet \mu_0 \circ \Phi_v
=0$ implying that $(A,\mu_0)$ is a $v$-algebra. We have, for any bilinear map $\varphi$ on $A$:
$$\delta^2_{v,\mu_0} \varphi (x,y,z)=\delta^2_{H,\mu_0} \varphi \circ \Phi_v(x,y,z) =\psi(xy,z)+\psi(yz,x)+\psi(zx,y)$$
for any $x,y,z \in A$, where $\psi$ is the skew-symmetric map associated to $\varphi$. Let $\mu_t=\mu_0+\sum t^i\varphi_i$ be a $v$-formal deformation of $\mu_0$. Since $v_{Lad} \in F_v$ (more precisely
$v_{Lad}\, v=3v_{Lad}$),
$\varphi_1$ is a Lie-admissible multiplication. Moreover $\delta^2_{v,\mu_0}\varphi_1=0$ and the Lie bracket $\psi_1$ satisfies
$$\psi_1(xy,z)+\psi_1(yz,x)+\psi_1(zx,y)=0$$
for any $x,y,z \in A.$ 
\begin{proposition}
Consider  the vector $v=Id+c+c^2$ of $\KS$ and  $(A,\mu_0)$  a commutative algebra. Then $(A,\mu_0)$ is a $v$-algebra and for any $v$-formal deformation $\mu=\mu_0+t\varphi_1+t^2\varphi_2 +\cdots$ of $\mu_0$, $(A,\psi_1)$ is a Lie  algebra such that
\begin{equation}\label{r2}
\psi_1(xy,z)+\psi_1(yz,x)+\psi_1(zx,y)=0
\end{equation}
for any $x,y,z \in A,$ where $\psi_1$ is the skew-symmetric bilinear map attached to $\varphi_1$.
\end{proposition}
 As the Leibniz rule implies Equation (\ref{r2}) we just have the conditions of a nonassociative Poisson algebra but  replacing the  Leibniz identity by a weak Leibniz identity (\ref{r2}) and we can define a notion of  $v$-Poisson and nonassociative $v$-Poisson algebras:
 
\begin{definition}
A nonassociative $v$-Poisson algebra is $\K$-vector space $A$ with  a Lie bracket $\psi$ and   a commutative multiplication $\mu$ tied up by the $v$-Leibniz identity:
$$\mathcal{L}(\mu,\psi)\circ \Phi_{v}=0.$$
A  $v$-Poisson algebra $(A,\mu,\psi)$ is  a  nonassociative $v$-Poisson  such that $\mu$  is moreover associative.
 \end{definition}
We trivially have that a Poisson algebra is a $v$-Poisson algebra and also a nonassociative
$v$-Poisson algebra.
We then obtain 
\begin{proposition}
Let $\mu=\mu_0+t\varphi_1+t^2\varphi_2 + \cdots$ be a $(Id+c+c^2)$-formal deformation of a commutative  multiplication 
$\mu_0$.Then
 $(A[[t]],\mu_t)$ is a deformation quantization of the nonassociative $(Id+c+c^2)$-Poisson algebra $(A,\mu_0,\psi_1)$. 
\end{proposition}
If $\mu_0$ is commutative associative, $(A,\mu_0, \psi_1)$ is a $(Id+c+c^2)$-Poisson algebra. Then
\begin{corollary}
Let $\mu=\mu_0+t\varphi_1+t^2\varphi_2 + \cdots$ be a $(Id+c+c^2)$-formal deformation of a commutative  associative multiplication 
$\mu_0$.Then
 $(A[[t]],\mu_t)$ is a deformation quantization of the $(Id+c+c^2)$-Poisson algebra $(A,\mu_0,\psi_1)$. 
\end{corollary}

 An example of algebra which is nonassociative $(Id+c+c^2)$-Poisson algebra but not nonassociative Poisson is obtained by considering the $2$-dimensional case: let $\{e_1,e_2\}$ be a basis of $A$ and
$$\psi_1(e_1,e_2)=e_2, \ \ \mu_0(e_1,e_1)=2\beta e_1, \ \mu_0(e_1,e_2)=\mu_0(e_2,e_1)=\alpha e_1+\beta e_2, \ \mu_0(e_2,e_2)=2\alpha e_2.$$
The algebra $(A,\mu_0,\psi_1)$ is a nonassociative $(Id+c+c^2)$-Poisson and it is a Poisson algebra when $\alpha=0$ and $\beta=1.$ Recall  (see \cite{RM}) that a Poisson algebra $(A, \mu, \{\, , \} )$ is also represented by only one multiplication $\cdot$ which satisfies a nonassociative identity 
$$x\cdot (y\cdot z)-(x\cdot y )\cdot z+\frac{1}{3} \left[ (x\cdot z )\cdot y+(y\cdot z )\cdot x-(y\cdot x )\cdot z-(z\cdot x )\cdot y\right] =0$$
and the two multiplications $\mu$ and $\{ \, , \, \}$ appearing in the definition of Poisson algebras are reobtained  by the depolarization process. This nonassociative multiplication is called Poisson admissible.  If we apply this idea to nonassociative $(Id+c+c^2)$-Poisson algebras, we find that the class of nonassociative  $(Id+c+c^2)$-Poisson admissible algebras corresponds to the 3-power associative algebras $(A,\cdot)$ that is the multiplication $\cdot$ satisfies 
$$\mathcal{A}_{\, \cdot} \circ \Phi_{v_{3Pa}}=0.$$
There is a one-to-one correspondence between nonassociative $(Id+c+c^2)$-Poisson algebras and 3-power associative algebras (see Section 7).

\subsection{Rank$(v)=3 : v=\alpha Id-\alpha \tau_{12} +(\alpha+\beta-3)\tau_{13}-\beta \tau_{23}+\beta c +(3-\alpha -\beta) c^2$ with $(\alpha,\beta) \neq (1,1)$ }

We will focus in this section on Vinberg, Pre-Lie and $G_3$-associative algebras after studying the general case which shows that the cases where $\beta=1$ and $\alpha \neq 1$ have additional properties on $\psi_1$, the skew-symmetric multiplication associated to $\varphi_1$ and  so are particular in this family.

\subsubsection{General case $v=\alpha Id-\alpha \tau_{12} +(\alpha+\beta-3)\tau_{13}-\beta \tau_{23}+\beta c +(3-\alpha -\beta) c^2$ with $(\alpha,\beta) \neq (1,1)$} As in previous cases, if $\mu_t$ is a  $v$-formal deformation of a  commutative $v$-associative multiplication  $\mu_0$, then $\varphi_1$ is Lie-admissible and its commutator $\psi_1$ is a Lie bracket. The equation $\delta_{v,\mu_0}^2\varphi_1=0$ gives additional properties on $\psi_1$ if only if $\beta=1$. In fact
$$
\begin{array}{rl}
    \delta_{v,\mu_0}^2\varphi_1(x,y,z)  =&x((\alpha-\beta) \varphi_1(y,z)-(\alpha+2\beta-3)\varphi_1(z,y)) + z(2\alpha+\beta-3)\psi_1(y,x)  \\
      & +y((\alpha+2\beta-3)\varphi_1(z,x)-(\alpha-\beta)\varphi_1(x,z)) +(\alpha-\beta)\varphi_1(x,yz)\\
     & -(\alpha+2\beta-3)\varphi_1(yz,x)+(\beta-\alpha)\varphi_1(y,xz)-(3-\alpha-2\beta)\varphi_1(zx,y)\\
   =  &(\alpha+2\beta-3)(x\psi_1(y,z)-z\psi_1(x,y)+y\psi_1(z,x)-\psi_1(yz,x)+\psi_1(zx,y))\\
   & -3(\beta-1)(x\varphi_1(y,z)-y\varphi_1(x,z)+\varphi_1(x,yz)-\varphi_1(y,xz))\\
   &-(\alpha-\beta)z\psi_1(x,y).
\end{array}
$$
Then $ \delta_{v,\mu_0}^2\varphi_1(x,y,z)  =0$ gives a relation concerning only $\psi_1$ as soon as $\beta=1$. If $\beta \neq 1$ we have to consider additional condition $x\varphi_1(y,z)-y\varphi_1(x,z)+\varphi_1(x,yz)-\varphi_1(y,xz))=0$ which doen't concern all cocycles $\varphi_1$. Then we assume $\beta=1$. Since we assumed that $(\alpha,  \beta) \neq (1,1)$,
 then $\alpha \neq 1$. Because this hypothesis,  the $G_i$-algebras for $i=2,3,4$ are excluded In fact $G_2$-algebras corresponds to $\alpha=3,\beta=0$, $G_3$-algebras to $\alpha=0,\beta=3$  and $G_4$-algebras to $\alpha=0,\beta=0$; we will see later on some relations on $\varphi_1$ or on $\rho_1$. For $v$-associative algebras with  $v=\alpha (Id-\tau_{12}) - \tau_{23}+ c +(2-\alpha ) (c^2-\tau_{13})$ and $\alpha \neq 1$,   the equation  $\delta_{v,\mu_0}^2\varphi_1=0$  reduces to
$$\psi_1(x,yz)-\psi_1(y,xz) +x\psi_1(y,z)-y\psi_1(x,z)-2z\psi_1(y,x)=0.$$
Considering the Leibniz operator
$$\mathcal{L}(\mu_0,\psi_1)(x,y,z)=\psi_1(xy,z)-x\psi_1(y,z)-\psi_1(x,z)y,$$
the equation  $\delta_{v,\mu_0}^2\varphi_1=0$ is then equivalent to
$$\mathcal{L}(\mu_0,\psi_1)\circ \Phi_{Id-\tau_{12}}=0.$$
We then have
\begin{proposition}
If $\mu_t$ is a $v$-deformation of a commutative $v$-associative algebra $(A,\mu_0)$ with $v=\alpha Id-\alpha \tau_{12} +(\alpha-2)\tau_{13}- \tau_{23}+ c +(2-\alpha) c^2$ and $\alpha \neq 1$, then $(A,\mu_0,\psi_1$) is a nonassociative  $(Id-\tau_{12})$-Poisson algebra and $(A[[t]],\mu_t)$ is a deformation quantization of this nonassociative  $(Id-\tau_{12})$-Poisson algebra.
\end{proposition}

\subsubsection{$G_2$-algebras or Vinberg algebras} For $v=Id-\tau_{12},$ a $v$-algebra is also called a Vinberg algebra. Let $\mu_0$ be a commutative Vinberg algebra. For any bilinear map $\varphi$, we have
$$\delta^2_{v,\mu_0}\varphi(x,y,z)=-x\varphi(y,z)+y\varphi(x,z)-\varphi(x,yz)+\varphi(y,xz)+z\psi(x,y)$$
for any $x,y,z \in A$ with $\psi$ the skew-symmetric bilinear map   attached to $\varphi$.
Let $\mu_t=\mu_0+\sum t^i\varphi_i$ be a $(Id-\tau_{12})$-formal deformation of $\mu_0$. Using the same notations as above, we have
$$\left\{
 \begin{array}{l}
   \delta^2_{v,\mu_0}\varphi_1=0,    \\
     \varphi_1 \bullet_v \varphi_1 +   \delta^2_{v,\mu_0}\varphi_2=0. 
\end{array}
\right.
$$
As $v_{Lad} \in F_v$, the multiplication $\varphi_1$ is Lie-admissible the algebra $(A,\psi_1)$ is a Lie algebra. 
% Now, look for the identities concerning the Lie bracket $\psi_1$. Since 
%$\delta^2_{v,\mu_0}\varphi_1=0$, we have $\delta^2_{v_1v,\mu_0}\varphi_1=0$ for any $v_1=\vv$. But 
%$$v_1v=(a_1-a_2)(Id-\tau_{12})+(a_3-a_5)(\tau_{13}-c)+(a_4-a_6)(\tau_{23}-c^2).$$
%implying that $v_1v \in F_v$ for any $v_1$. 
% Then the relation $\delta^2_{v_1v,\mu_0}\varphi_1=0$ is equivalent to $\delta^2_{v,\mu_0}\varphi_1=0$ and doen't give any new information. Nevertheless, we can write again the expression

 The equation $\delta^2_{v,\mu_0}\varphi_1=0$ writes :
%$$\mathcal{L}(\mu_0,\varphi_1)(x,y,z)=\varphi_1(\mu_0(x,y),z)-\mu_0(x,\varphi_1(y,z))-\mu_0(\varphi_1(x,z),y).$$ We obtain
$$\delta^2_{v,\mu_0}\varphi_1=-\mathcal{L}_R(\mu_0,\varphi_1) \circ \Phi_{Id-\tau_{12}}=0,$$
using the right-Leibniz operator:
 $$\mathcal{L}_R(\mu_0,\varphi_1)(x,y,z)=\varphi_1(x,\mu_0(y,z))-\mu_0(y,\varphi_1(x,z))-\mu_0(\varphi_1(x,y),z).$$
\begin{proposition}
Let $(A,\mu_0)$ be a commutative Vinberg algebra. Then any  $(Id-\tau_{12})$-formal deformation $\mu=\mu_0+t\varphi_1 + \cdots$ determines an Lie-admissible  algebra $(A,\varphi_1)$ satisfying
$$\mathcal{L}_R(\mu_0,\varphi_1) \circ \Phi_{Id-\tau_{12}}=0,$$
\end{proposition}

\subsubsection{$G_4$-algebra also called Pre-Lie algebras: $v=Id-\tau_{23}$}
This case is similar to the $G_2$-algebra case: 
\begin{proposition}
Let $(A,\mu_0)$ be a commutative Pre-Lie algebra. Then any formal $(Id-\tau_{23})$-formal deformation $\mu=\mu_0+t\varphi_1 + \cdots$ determines a Lie-admissible  algebra $(A,\varphi_1)$ satisfying
$$\mathcal{L}(\mu_0,\varphi_1) \circ \Phi_{Id-\tau_{23}}=0.$$

\end{proposition}

\subsubsection{$G_3$-algebra: $v=Id-\tau_{13}$}
A commutative $v$-algebra is also associative. As $v_{Lad} \in F_v$, the linear term $\varphi_1$ of a $v$-formal deformation of $\mu_0$ is Lie-admissible. The map $\varphi_1$ satisfies also
$$\delta^2_{v,\mu_0}\varphi_1(x,y,z)=-x\rho_1(y,z)+z\rho_1(x,y)+\rho_1(xy,z)-\rho_1(x,yz)=0$$where $\rho_1$ is the symmetric map attached to $\varphi_1$, which can also be written 
$$\delta^2_{v,\mu_0}\varphi_1(x,y,z)=-\rho_1(x,yz)+z\rho_1(x,y)+y\rho_1(x,z)+\rho_1(xy,z)-x\rho_1(y,z)-y\rho_1(x,z)=0,$$
that is
$$\LE(\mu_0,\rho_1)(x,y,z)-\LE(\mu_0,\rho_1)(z,y,x)=0.$$
because
$$\LE(\mu_0,\rho_1)(x,y,z)=\LE_R(\mu_0,\rho_1)(z,y,x)$$
as $\rho_1$ is a commutative multiplication.

\noindent Remark that il $\mu_0$ is commutative and 
$\psi_1$ is skew symmetric, we have that 
$$\LE(\mu_0,\psi_1)(x,y,z)=-\LE_R(\mu_0,\psi_1)(z,y,x).$$

\begin{proposition}

Let $(A,\mu_0)$ be a commutative  $(Id-\tau_{13})$-algebra. Then, if $\mu=\mu_0+t\varphi_1+t^2\varphi_2+\cdots$ is a $(Id-\tau_{13})$-formal deformation of $\mu_0$, then if $\psi_1$ and $\rho_1$ are respectively the skew-symmetric and symmetric bilinear maps associated to $\varphi_1$
\begin{enumerate}
\item $(A,\psi_1)$ is a Lie algebra,
\item The symmetric map $\rho_1$ satisfies
$$\LE(\mu_0,\rho_1)\circ \Phi_{Id-\tau_{13}}=0.$$
\end{enumerate}
\end{proposition}

\medskip

%\noindent{Example.} Let $A=\mathcal{F}^1(\R,\R)$ be the commutative algebra of derivable real functions. Since it is associative, it is also $(Id-\tau_{13})$-associative for the natural product $\mu_0$. Let $\varphi_1$ the bilinear ap on $A$:
%$$\varphi_1(f,g)=fg'.$$
%Then $\psi_1(f,g)=fg'-f'g$ is a Lie braket and $\rho_1(f,g)=f'g+fg'=(fg)'$.  We have also
%$$f\rho_1(g,h)-h\rho_1(f,g)-\rho_1(fg,h)+\rho_1(f,gh)=f(gh)'-h(fg)'-fgh'+f(gh)'=0.$$

\subsection{Rank$(v)=4$: $v=2Id+(1+\alpha) \tau_{12} + \tau_{13} +c+(1-\alpha)c^2$ with $\alpha \neq 1$}
Let $\mu_t=\mu_0+t\varphi_1+ \cdots$ be a $v$-formal deformation of a commutative $v$-associative algebra. Then $(A,\varphi_1)$ is a Lie-admissible algebra and the equation $\delta^2_{\mu_0,v}\varphi_1=0$ is equivalent to a quadratic relation on 
%$\varphi_1$ or 
$\psi_1$ if only if $\alpha =-\frac{1}{2}$ or $\delta^2_{H,\mu_0}\varphi_1=0.$ In fact 
$$
\begin{array}{ll}
   \delta^2_{v,\mu_0}\varphi_1(x,y,z)=   &  x\psi_1(y,z)-\alpha\ y\psi_1(z,x) +z \psi_1(x,y)  +(2-\alpha)\psi_1(z,xy)-2\psi_1(yz,x)\\
      &   +(-1+\alpha)\psi_1(zx,y)+(1+2\alpha)\delta^2_{H,\mu_0}\varphi_1(y,x,z).\\
\end{array}
$$
Then $ \delta^2_{v,\mu_0}\varphi_1=0$ implies a quadratic relation on $\psi_1$ as soon as  $(1+2\alpha)\delta^2_{H,\mu_0}\varphi_1=0$. Since $\delta^2_{H,\mu_0}\varphi_1=0$ is a particular case of   $\delta^2_{v,\mu_0}\varphi_1=0$, and since we want a generic identity, then $\alpha =-\frac{1}{2}$. This corresponds to 
weakly-associative algebra also called Lie-admissible flexible algebras \cite{RWass} with vector $v=Id-\tau_{12}+c$. We recall results obtained in \cite{RWass}. 
Let  $\mu_0$ be a commutative  associative algebra.  Then it is also $v$-associative and we can consider a $v$-formal deformation of $\mu_0$:
$$\mu=\mu_0+t\varphi_1 + t^2 \varphi_2 + \cdots$$
We deduce
\begin{equation}\label{eqv}
\left\{
\begin{array}{ll}
\medskip
 {\rm order} \  0 & \mathcal{A}_{\mu_0} \circ \Phi_v =0,      \\
 \medskip
   {\rm order} \  1   &   \delta_{v,\mu_0 }^2\varphi_1 =\delta_{H,\mu_0 }^2\varphi_1 \circ \Phi_v=0,\\
   \medskip
   {\rm order} \  2  & \varphi_1 \bullet_v \varphi_1 +  \delta_{H,\mu_0}^2 \varphi_2 \circ \Phi_v=0.
\end{array}
\right.
\end{equation}
%We therefore want to determine all vectors $v\in \KS$ such that any $v$-deformation determines a Poisson algebra $(A,\mu_0,\psi_1)$. 
%The skew bilinear map $\psi_1$ attached to $\varphi_1$ is  a Lie bracket, if and only if  $\varphi_1$ is Lie-admissible that is $ \varphi_1 \bullet_{v_{Lad}} \varphi_1=  \varphi_1 \bullet \varphi_1 \circ \Phi_{v_{Lad}}=0$.  But Equations (\ref{eqv}) imply that $\varphi_1 \bullet_v \varphi_1  =- \delta_{H,\mu_0}^2 \varphi_2 \circ \Phi_v$. From Proposition \ref{Mv}, there exists  $v_1 \in \KS$ be such that 
%$$
%\left\{
%\begin{array}{l}
%v_1v=v_{Lad},\\
% \delta_{H,\mu_0}^2 \varphi_2 \circ \Phi_{v_{Lad}}=\delta_{H,\mu_0}^2 \varphi_2 \circ \Phi_{v_1v} =0,\\
% \end{array}
% \right.
% $$
%if and only if $v_{Lad} \in F_v$. 
%Such vectors $v$ are described in \cite{G.R.Nonass}. 
From Proposition \ref{Mv}, $\varphi_1$ is Lie-admissible and $\psi_1$ is a Lie bracket.
Let us  now investigate the consequences of the equation $\delta_{H,\mu_0 }^2\varphi_1 \circ \Phi_v=0$ by considering a vector $w$ canceling the $\rho_1$, that is $w=a_1Id+a_2\tau_{12}+a_3\tau_{13}+a_4\tau_{23}+(a_1-a_3+a_4)c+(a_1+a_2-a_3)c^2.$  We have using Leibniz operator
 associated to $\mu_0$ and $\psi_1$   that $\delta_{H }^2\varphi_1\circ \Phi_w =0$ is equivalent to
 $$
\begin{array}{r}
(a_1+a_2)\LE(\mu_0,\psi_1) (x,y,z)  +(a_1+a_4)\LE(\mu_0,\psi_1) (y,z,x) \\
  +(a_1+a_2-a_3+a_4)\LE(\mu_0,\psi_1) (z,x,y)=0
\end{array}
$$
If the components of $w$ satisfies one of the following conditions
\begin{enumerate}
  \item $a_1+a_2 \neq 0$ and $a_1+a_4=a_2-a_3=0$,
  \item  or $a_1+a_4 \neq 0$ and $a_1+a_2=a_3-a_4=0$,
   \item or $a_3-a_4 \neq 0$ and $a_1+a_2=a_1+a_4=0.$
\end{enumerate}
then $\delta_{H }^2\varphi_1\circ \Phi_w=0$ implies $\LE(\mu_0,\psi_1)=0.$ 
For each case, the $w$ vector belongs to  $F_{Id-\tau_{12}+c}.$ In fact, let us consider the first case.  The equation $(\alpha_1Id+\alpha_2\tau_{12}+\alpha_3\tau_{13}+\alpha_4\tau_{23}+\alpha_5 c+\alpha_6 c^2)\circ (Id-\tau_{12}+c)=w$ is equivalent to the linear system
$$
\left\{
\begin{array}{l}
    \alpha_1-\alpha_2+\alpha_6=a_1   \\
       -\alpha_1+\alpha_2+\alpha_3=a_2   \\  
          \alpha_3+\alpha_4-\alpha_5=a_2   \\
             \alpha_2+\alpha_4-\alpha_6=-a_1   \\
                \alpha_1-\alpha_3+\alpha_5=-a_2   \\
                  - \alpha_4+\alpha_5+\alpha_6=a_1   \\
\end{array}
\right.
\Leftrightarrow
\left\{
\begin{array}{l}
    \alpha_4=-\alpha_1   \\
       \alpha_5=-\alpha_2  \\  
          \alpha_6=a_1+\alpha_2-\alpha_1   \\
             \alpha_3=a_2-\alpha_2+\alpha_1 
\end{array}
\right.
$$
which has a nontrivial solution. It is similarly for the other cases. 
%Moreover $v_{Lad}$ belongs to $F_{Id-\tau_{12}+c}.$
Then we can find a vector $v_1$ such that $v_1v=w$ and $\delta_{H }^2\varphi_1\circ \Phi_w=0$ implies $\LE(\mu_0,\psi_1)=0.$ This implies:

\begin{proposition} Let $(A,\mu_0)$ be a commutative algebra and consider the vector $v=Id-\tau_{12}+c $. Then $(A,\mu_0)$ is $v$-associative and for any $v$-formal deformation $\mu_t=\mu_0+t\varphi_1+t^2\varphi_2 + \cdots$ of $\mu_0$, the algebra $(A,\mu_0,\psi_1)$ is a nonassociative Poisson algebra where $\psi_1$ is the skew-symmetric application associated with $\varphi_1$.
\end{proposition}
\begin{corollary} \label{coro 18}
Let $(A,\mu_0)$ be an associative commutative algebra and consider the vector $v=Id-\tau_{12}+c $. Then $(A,\mu_0)$ is $v$-associative and for any $v$-formal deformation $\mu_t=\mu_0+t\varphi_1+t^2\varphi_2 + \cdots$ of $\mu_0$, the algebra $(A,\mu_0,\psi_1)$ is a Poisson algebra.
\end{corollary}

\noindent{\bf Consequence.} In the usual deformation quantization process, a Poisson algebra $(A,\mu_0,\psi_1)$  is obtained from a  formal deformation $\mu_t=\mu_0+t\varphi_1+ \cdots $ of a commutative associative  algebra and the algebra $(A[[t]],\mu_t)$ is a deformation quantization of the Poisson algebra $(A,\mu_0,\psi_1)$. Corollary \ref{coro 18} shows that a Poisson algebra  $(A,\mu_0,\psi_1)$ is also obtained from a $(Id-\tau_{12}+c)$-formal deformation  (that is a weakly associative formal deformation) $\mu_t=\mu_0+t\varphi_1+ \cdots $ of a commutative associative  algebra. Thus  we also consider  the algebra $(A[[t]],\mu_t)$  as a quantization of a Poisson algebra  $(A,\mu_0,\psi_1)$ in this more general case. The $v$-algebras with $v=Id-\tau_{12}+c$ called weakly-associative algebras have been introduced  in \cite{RWass} where an algebraic study is presented.

\medskip

\noindent{\bf Remark.} In \cite{RWass} we show that any commutative algebra, any Lie algebra, any associative algebra is weakly associative. In fact the class of weakly associative algebras (associated to $v=Id+c-\tau_{12}$) is the biggest  class containing the Lie algebras and the associative algebras such that the $v$-deformation of a commutative associative algebra gives a Poisson algebra so quantizations of a Poisson algebra. In \cite{RWass2}, we show also that the symmetric Leibniz algebras are also weakly associative.

\subsection{Rank$(v)=5$: $F_v=F_{2Id-\tau_{12}-\tau_{13}-\tau_{23}+c}$}
Let $\mu_t=\mu_0+\sum t^i \varphi_i$ be a $v$-formal deformation of the commutative $v$-algebra $(A,\mu_0)$. The product $\varphi_1$ is Lie-admissible  from Theorem \ref{co} and $\psi_1$ is a Lie bracket. Since $\mu_0$ is commutative, we have
$$\begin{array}{ll}
\delta^2_{v,\mu_0}\varphi_1(x,y,z)=& -x\varphi_1(y,z)-y\varphi_1(z,x)+2z\varphi_1(x,y)\\
 & +\varphi_1(xy,z)-\varphi_1(xz,y)-\varphi_1(x,yz)+\varphi_1(z,xy).
\end{array}
$$
Choosing $v_1=(a_1,a_2,a_3,a_4,a_1+a_2-a_3,a_1+a_2-a_4)$, we see that $\delta^2_{v,\mu_0}\varphi_1\circ \Phi_{v_1}$ contains only elements in $\psi_1$ so that $\delta^2_{v,\mu_0}\varphi_1\circ \Phi_{v_1}=0$ reads
$$
\begin{array}{l}
(a_2-2a_3+a_4)x\psi_1(z,y)+(-a_2-a_3+2a_4)y\psi_1(z,x)+(2a_2-a_3-a_4)z\psi_1(x,y) \\
+(a_4-a_3)\psi_1(xy,z)+(a_3-a_2)\psi_1(xz,y)+(a_2-a_4)\psi_1(yz,x)=0
\end{array}
$$
that is 
$$
\begin{array}{l}
a_2(-x\psi_1(y,z)-y\psi_1(z,x)+2z\psi_1(x,y) -\psi_1(xz,y)+\psi_1(yz,x))\\
+a_3(-2x\psi_1(y,z)-y\psi_1(z,x)-z\psi_1(x,y)-\psi_1(xy,z)+\psi_1(xz,y)\\
+a_4(x\psi_1(y,z)+2y\psi_1(z,x)-z\psi_1(x,y)+\psi_1(xy,z)-\psi_1(yz,x)=0.
\end{array}
$$
This is equivalent to the identity
$$-x\psi_1(y,z)-y\psi_1(z,x)+2z\psi_1(x,y) -\psi_1(xz,y)+\psi_1(yz,x)=0$$
that is 
$$\LE(\mu_0,\psi_1)(y,z,x)-\LE(\mu_0,\psi_1)(x,z,y)=0.$$
\begin{proposition}
Let $(A,\mu_0)$ be a commutative $v$-algebra with $v=2Id-\tau_{12}-\tau_{13}-\tau_{23}+c.$ Any $v$-formal deformation $\mu_t=\mu_0+\sum t^i \varphi_i$ of $(A,\mu_0)$ is a deformation quantization of a nonassociative $(Id-\tau_{13})$-Poisson algebra $(A,\mu_0,\psi_1)$ that is
\begin{enumerate}
\item $\mu_0$ is a commutative multiplication on $A$,
\item $\psi_1$ is a Lie bracket on $A$,
\item $\LE(\mu_0,\psi_1)\circ \Phi_{Id-\tau_{13}}=0.$
\end{enumerate}
\end{proposition} 

%%%%%%%%%%%%%%%%%%
\section{A generalization : $(\KS)^2$-associative algebras}
This generalization has been introduced in \cite{G.R.Nonass}. If $(A,\mu)$ is a $\K$-algebra, we denote by $\mathcal{A}_{\mu}$ the associator of $\mu$. Let’s write this associator in the following form:
$$ \mathcal{A}_{\mu}=\mathcal{A}_{\mu}^L-\mathcal{A}_{\mu}^R
$$
where $\mathcal{A}_{\mu}^L=\mu\circ (\mu \otimes Id)$ and $\mathcal{A}_{\mu}^R=\mu \circ(Id \otimes \mu).$

Now, instead of considering action of $\Sigma_3$-permutation on the associator we can consider it independently on $\mathcal{A}_{\mu}^L$ and $\mathcal{A}_{\mu}^R$ which will induce different symmetries.

\subsection{Definition}

\begin{definition}
Let $v$ and $w$ be two vectors of $\KS$. We say that the algebra $(A,\mu)$ is a $(v,w)$-algebra if we have
$$(\star)
%\left\{
%\begin{array}{l}
%\mathcal{A}_{\mu}^L \circ \Phi_v=0 \\
%\mathcal{A}_{\mu}^R \circ \Phi_w=0
%\end{array}
%\right.
%\quad {\rm or} \quad (\star\star)
\left\{
\begin{array}{l}
\mathcal{A}_{\mu}^L \circ \Phi_v - \mathcal{A}_{\mu}^R \circ \Phi_w=0
\end{array}
\right.
$$
If there exists a non trivial  $(v,w) \in \KS^2 $  such that $(A,\mu)$ is a $(v,w)$-algebra, the algebra $(A,\mu)$ is $(\KS)^2$-associative.
\end{definition}

%We assume also that for every $v',w' \in \KS$
%such that $v \in F_{v'}, \ w \in F_{w'}$ and $v' \notin F_v,\ w' \notin F_w$ we have $\mathcal{A}_{\mu}^L \circ \Phi_{v'} \neq 0$
%or $\mathcal{A}_{\mu}^R \circ \Phi_{w'}\neq 0.$ 

An interesting case of $(v,w)$-algebras corresponds to the algebras given by the system of equations
$$(\star\star)
\left\{
\begin{array}{l}
\mathcal{A}_{\mu}^L \circ \Phi_v=0, \\
\mathcal{A}_{\mu}^R \circ \Phi_w=0.
\end{array}
\right.
$$

In the study of $(v,w)$-algebras, the first study is to know if a $(v,w)$-algebra can be defined as a $v_1$-algebra. The easiest example is when $v=w$. In \cite{G.R.Nonass} we study the $(v,w)$-algebras which are Lie-admissible algebras or Pre-Lie algebras. 

\medskip

\noindent{Example.} A Leibniz algebra satisfies the quadratic relation
$$
\mu(\mu(x,y),z)=\mu(x,\mu(y,z))+\mu(\mu(x,z),y)$$
 or
$$\mathcal{A}_{\mu}^L \circ \Phi_{Id-\tau_{23}}-\mathcal{A}_{\mu}^R\circ \Phi_{Id}= 0.$$
Then Leibniz algebras are $(Id-\tau_{23},Id)$-algebras. 
Symmetric Leibniz algebras are defined by a pair of quadratic relations. They correspond to
$$ \mathcal{A}_{\mu}^L\circ \Phi_{Id-\tau_{23}}-\mathcal{A}_{\mu}^R \circ\Phi_{Id} , \ \ \ \mathcal{A}_{\mu}^L \circ\Phi_{Id} -\mathcal{A}_{\mu}^R\circ \Phi_{Id-\tau_{12}},$$

The notion of $(v,w)$-formal deformation of a $(v,w)$-algebra is similar to this notion for $v$-algebras.  Let $(A,\mu_0)$ be a $(v,w)$-algebra defined by a relation of type ($\star,\star \star$). Consider  $\ds \mu_t= \mu_0 +\sum_{ i \geq1} t^i\varphi_i$.  We say that $\mu_t$ is a $(v,w)$-formal deformation of $\mu_0$ if $(A[[t]],\mu_t)$ is a $(v,w)$-algebra. To describe the relations between the $\varphi_i$, we need to introduce some notations:
$$\varphi_i \bullet_v^L \varphi _j = \varphi_i \circ (\varphi _j \otimes Id ) \circ \Phi_v, \ \ \varphi_i \bullet_w^R \varphi _j = \varphi_i \circ (Id \otimes \varphi _j) \circ \Phi_w$$
$$\delta^{2,L}_{v,\mu_0}\varphi=\varphi \bullet^L_v \mu_0+ \mu_0 \bullet^L_{v} \varphi, \ \ \delta^{2,R}_{w,\mu_0}\varphi=\varphi \bullet^R_w \mu_0+ \mu_0 \bullet^R_{w} \varphi.$$
Thus, to say that $\mu_t$ is a $(v,w)$-formal deformation of $\mu_0$ implies in particular:
\begin{enumerate}
\item order $0$ : $\mu_0$ is a $(v,w)$-algebra,
\item order $1$ : $\delta^{2,L}_{v,\mu_0}\varphi_1- \delta^{2,R}_{w,\mu_0}\varphi_1 =0$,
\item order $2$ : $\varphi_1 \bullet^L_v \varphi_1-\varphi_1 \bullet^R_w \varphi_1 + \delta^{2,L}_{v,\mu_0}\varphi_2
- \delta^{2,R}_{w,\mu_0}\varphi_2=0$
\end{enumerate}

\subsection{A fundamental example: the anti-associative algebras}
\begin{definition}
A $\K$-algebra $(A,\mu)$ is called anti-associative if the multiplication $\mu$ satisfies the following identity
$$\mu(\mu(x,y),z)+\mu(x,\mu(y,z))=0$$
for any $x,y,z \in A$.
\end{definition}
We will denote by $\mathcal{AA}_\mu$ the trilinear map
$$\mathcal{AA}_\mu (x,y,z)=\mu(\mu(x,y),z)+\mu(x,\mu(y,z))$$
and $(A,\mu)$ is \ant if and only if $\mathcal{AA}_\mu =0.$ In terms of $(v,w)$-algebra an \ant algebra is an $(Id,-Id)$-algebra.
 In an \ant algebra all the  $4$-products are zero. In fact
$$
\begin{array}{rl}
x(y(zt)) &=-x((yz)t)=(x(yz))t=-((xy)z)t   \\
{\rm and}  \, \,    x(y(zt))  &   =-(xy)(zt)=((xy)z)t
\end{array}
$$
so all $4$-products are trivial. Anti-associative
algebras are therefore always 3-step nilpotent. 

\noindent {\bf Examples}
There are `natural' examples of the
anti-associativity. For instance, the standard basis elements $\{1,e_1,e_2,e_3,e_4,e_5,e_6,e_7\}$
 of the octonions (also called the  Cayley algebra)
satisfy
\[
(e_ie_j)e_k=-e_i(e_je_k),
\]
whenever $e_i e_j \neq \pm e_k$ and $1 \leq i,j,k\leq 7$ are  distinct.
In \cite{M-R-galgalim} anti-associative algebras in small dimension are described. For example in dimension $3$ we have obtained the following non isomorphic nontrivial anti-associative algebras $(A,\cdot)$:
\begin{enumerate}
\item $e_i \cdot e_i=0, \ e_1 \cdot e_2=-e_2 \cdot e_1=e_3$
\item $e_1 \cdot e_1=e_2, \ e_1 \cdot e_2=-e_2 \cdot e_1=e_3$ which happens to be the free anti-associative algebra on one
generator,
\item $e_1 \cdot e_1=e_2,\ e_1 \cdot e_3=ae_2, \ e_3 \cdot e_1=be_2 ,\ e_3 \cdot e_3=e_2, $
\item $e_1 \cdot e_1=e_2,\ e_1 \cdot e_3=ae_2, \ e_3 \cdot e_1=be_2.$
\end{enumerate}
with $a,b \in \K$ and where $\{e_1,e_2,e_3\}$ is a basis of $A$.

\medskip

If $\mathcal{AA}ss$ denotes the quadratic operad corresponding to the \ant algebra, then
$$\dim \mathcal{AA}ss(1)=1, \ \dim \mathcal{AA}ss(2)=2, \ \dim \mathcal{AA}ss (3)= 6 \ {\rm and} \ \dim \mathcal{AA}ss (n)=0  \ {\rm for} \ n \geq 4.$$
Let us note also that this operad is self dual. In \cite{M-R-galgalim} it is proved that the operad  $\mathcal{AA}ss$ is not Koszul computing the inverse series of the generating function $g_{\mathcal{AA}ss}(t)=t+t^2+t^3$.

Concerning the problem of deformation of anti-associative algebras, the `standard' cohomology  of an
anti-associative algebra $A$ with coefficients in itself is described in \cite{M-R-galgalim}  and compared
to the relevant part of the deformation cohomology  based on
the minimal model of the anti-associative operad $\mathcal{AA}ss$. Since
$\mathcal{AA}ss$ is not Koszul, these two cohomologies
differ. The standard cohomology $H^*_{st}(A,A)$ is the cohomology of the complex
$$
C^1(A,A) \stackrel{\delta^1_{AA}}{\longrightarrow} C^2(A,A)
\stackrel{\delta^2_{AA}}{\longrightarrow}
C^3(A,A) \stackrel{\delta^3_{AA}}{\longrightarrow}
0 \stackrel{0}{\longrightarrow} 0 \stackrel{0}{\longrightarrow} \cdots
$$
in which $C^p(A,A) := Hom(A^{\otimes p},A)$ for $p = 1,2,3$, and all higher
$C^p$'s are trivial. The two nontrivial pieces of the differential are
basically the Hochschild differentials with ``wrong'' signs of some
terms:
\begin{align*}
\delta^1_{AA}(f)(x,y)& := x f(y) - f(xy) + f(x)y,
\mbox { and}
\\
\delta^2_{AA}(\varphi)(x,y,z) &:= x \varphi(y,z)+ \varphi(xy,z)+ \varphi(x,yz) +\varphi(x,y)z,
\end{align*}
for $f\in Hom(V,V)$, $\varphi \in Hom(A^{\otimes 2},V)$ and $x,y,z \in
V$. One sees, in particular, that 
$$H^*_{st}(A,A)^p = 0$$ for $p\geq 4$.

\medskip

The  deformation cohomology of anti-associative algebras, based of the study of a minimal model is also studied in \cite{M-R-galgalim}. We summarize the results:
we consider the complex 
\[
C^1_{AAss}(A,A) \stackrel{\delta^1}{\longrightarrow} C^2_{AAss}(A,A) 
\stackrel{\delta^2}{\longrightarrow}
C^3_{AAss}(A,A)  \stackrel{\delta^3}{\longrightarrow}
C^4_{AAss}(A,A)  \stackrel{\delta^4}{\longrightarrow} \cdots
\]

-- $C^1_{AAss}(A,A) = Hom(A,A)$

-- $C^2_{AAss}(A,A) = Hom(A^{\otimes 2},A)$

-- $C^3_{AAss}(A,A) = Hom(A^{\otimes 3},A)$, and

-- $C^4_{AAss}(A,A) = Hom(A^{\otimes 5},A) \oplus Hom(A^{\otimes 5},A)  \oplus
 Hom(A^{\otimes 5},A) \oplus Hom(A^{\otimes 5},A) $.

\noindent
Observe that $C^p_{AAss}(A,A)= C^p(A,A) $ for $p = 1,2,3$, while 
$C^4_{AAss}(A,A)$ consists of $5$-linear maps. The
differential $\delta^p$ agrees with $\delta_{AA}^p$ for  $p = 1,2$ while,
for $g \in C^3_{AAss}(A,A)$, one has
\[
\delta^3 (g) = (\delta_1^3(g),\delta_2^3(g),\delta_3^3(g),\delta_4^3(g)),
\]
where
\begin{align*}
\delta_1^3(g)(x,y,z,t,u) &:= 
xg(y,z,tu) - g(x,y,z(tu)) + (xy)g(z,t,u) - g(xy,zt,u) 
\\
&\
+ g(xy,z,t)u - g((xy)z,t,u) + g(x,y,z)(tu) - g(x,yz,tu),
\\
\delta_2^3(g)(x,y,z,t,u) &:= 
g((xy)z,t,u) - g(xy,z,t)u + g(x,y,zt)u - g(x,y(zt),u)
\\
&\
+ xg(y,zt,u) - g(x,y,(zt)u) + (xy)g(z,t,u) - g(xy,z,tu),
\\
\delta_3^3(g)(x,y,z,t,u) &:= 
g(x,yz,tu) - xg(yz,t,u) + g(x,(yz)t,u) - x(g(y,z,t)u)
\\
&\
+ g(x,y,zt)u - g(xy,z,t)u + (g(x,y,z)t)u - g(x(yz),t,u), \mbox { and}
\\
\delta_4^3(g)(x,y,z,t,u) &:= 
g(xy,zt,u) - g(x,y,(zt)u) + x g(y,zt,u) - g(x,y(zt),u)
\\
&\
+ (xg(y,z,t))u - g(x,yz,t)u + (g(x,y,z)t)u - g(xy,z,t)u,
\end{align*}
for $x,y,z,t,u \in V$.

\medskip

We consider now formal deformation of anti-associative algebras. 
If $\mu_t=\mu_0+ \sum t^i \varphi_i$ is an \ant formal deformation of $\mu_0$, then we have, denoting by $xy$ the product $\mu(x,y)$:
\begin{enumerate}
  \item in degree $0$ : $\mu_0$ is \ant \!\!,
  \item In degree $1$: $\varphi_1(x,yz)+x\varphi_1(y,z)+\varphi_1 (x,y)z+\varphi_1(xy,z)=0$
  \item In degree $2$: $$\varphi_2(x,yz)+x\varphi_2(y,z)+\varphi_2 (x,y)z+\varphi_2(xy,z)+\varphi_1(x,\varphi_1(y,z))+\varphi_1(\varphi_1 (x,y),z) =0.$$
\end{enumerate}
that is
$$\delta_{AA,\mu_0}^2(\varphi_1)=0$$
and 
$$\delta_{AA,\mu_0}^2(\varphi_2)+\mathcal{AA}_{\varphi_1}=0. $$
Assume moreover that $\mu_0$ is commutative. 
It is not difficult to see that $\delta_{AA,\mu_0}^2(\varphi_2) \circ \Phi_v =0$ implies $v=0$ and we 
do not have a good deformation quantization frame work similar to the associative or $v$-associative cases. 
%have no algebras  such that formal deformations of $\mu_0$ are quantizations of this algebra.

\medskip

We are therefore naturally led to consider formal  deformation of a \ant product $\mu_0$ which is also skew-symmetric. In this case, we get an anti-commutative version of Lemma \ref{1} which was for  a commutative product.

\begin{lemma}\label{2}
Let $(A,\mu_0)$ be a anti-commutative algebra with $\mu_0 \neq 0$ and $\delta^2_{H,\mu_0}$  the coboundary operator:
$$\delta^2_{AA,\mu_0} \varphi (X,Y,Z)=X\varphi(Y,Z)+\varphi(XY,Z)+\varphi(X,YZ)+\varphi(X,Y)Z$$
where $X,Y,Z\in A$, the map $\varphi$ is bilinear on $A$  and  $XY$denotes the product $\mu_0(X,Y)$.  
Then $ \delta^2_{AA,\mu_0}\varphi \circ \Phi_{v_{3Pa}}= 0.$
\end{lemma}

\pf A direct computation proves, like for Lemma \ref{1} for the commutative  case that  $ \delta^2_{AA,\mu_0}\varphi \circ \Phi_{v_{3Pa}}= 0.$ \quad $\square$

%Let us note that $\varphi$, $\delta_{AA,\mu_0}^2(\varphi)\circ \Phi_{v}=0$
%implies $v=\alpha(1,1,1,1,1,1)$ that is $v=\alpha v_{3Pa}$. 
 \noindent As a consequence,$$\mathcal{AA}_{\varphi_1}\circ \Phi_{v_{3Pa}}=0.$$ 
We can say that $\varphi_1$ verify the  \ant  version of the $3$-power associative property. We deduce

\begin{proposition}
Let $(A,\mu_0)$ be a skew-symmetric \ant algebra and $\mu= \mu_0+ \sum t^k \varphi_k$ an \ant 
formal deformation of $\mu_0$. Then,  if  $\rho_1$ denotes the symmetric part of $\varphi_1$, the algebra  $(A,\rho_1)$ is a Jacobi-Jordan algebra, that is 
\begin{enumerate}
\item $(A,\rho_1)$ is a commutative algebra,
\item $\rho_1$ satisfies the "Jacobi" identity:
$$\rho_1(x,\rho_1(y,z))+\rho_1(y,\rho_1(z,x))+\rho_1(z,\rho_1(x,y))=0$$
for all $x,y,z \in A$.
\end{enumerate}
\end{proposition}
\pf The Jacobi identity for $\rho_1$ follows from the fact that $\mathcal{AA}_{\varphi_1}\circ \Phi_{v_{3Pa}}=0.$

Let us examine the first condition $\delta_{AA,\mu_0}^2(\varphi_1)=0.$ With a similar proof to the associative case, we show that this identity implies
$$\rho_1(xy,z)+x\rho_1(y,z)+\rho_1(x,z)y=0$$
where $\mu_0(x,y)=xy.$
\begin{theorem}\label{antiP}
Let $(A,\mu_0)$ be a skew-symmetric \ant algebra and $\mu_t= \mu_0+ \sum t^k \varphi_k$ an \ant formal deformation of $\mu_0$. Then $(A,\mu_0,\rho_1)$ is an anti-Poisson algebra, that is
\begin{enumerate}
\item $(A,\rho_1)$ is a Jacobi-Jordan algebra,
\item The products $\mu_0$ and $\rho_1$ are tied up by the graded Leibniz identity:
$$\mathcal{L}_g(\mu_0,\rho_1)(x,y,z)=\rho_1(xy,z)+x\rho_1(y,z)+\rho_1(x,z)y=0.$$
\end{enumerate}
We will say that $\mu_t$ is a deformation quantization of the anti-Poisson algebra $(A,\mu_0,\rho_1)$.
\end{theorem}
Recall that an antiderivation of an algebra $(A,\mu)$ is a linear map $f$ such that
$$f(\mu(x,y))+\mu(x,f(y))+\mu(f(x),y)=0$$
for any $x,y$ in $A$. The graded Leibniz identity can be interpreted saying that for any $z \in A$, the linear maps $x \rightarrow \rho(x,z)$ is an antiderivation of the algebra $(A,\mu_0).$

%%%%%%%%%%%%%%

\noindent{\bf Remarks.} Considering an anti-commutative product $\mu_0$  and $w=a_1Id+a_2\tau_{12}
+a_3\tau_{13}+a_4\tau_{23}+a_5c+a_6c^2$ be a vector of $\KS$. The anti-commutativity of  $\mu_0$ implies \begin{equation}
\label{pp}
\begin{array}{l}
\delta^2_{AA,\mu_0}\varphi \circ \Phi_w(X,Y,Z)= \\

(a_1-a_5)X\varphi(Y,Z) +(a_4-a_3)X\varphi(Z,Y)+(a_2-a_4)Y\varphi(X,Z)+(a_5-a_6)Y\varphi(Z,X)\\
+(a_6-a_1)Z\varphi(X,Y)+(a_3-a_2)Z\varphi(Y,X)+(a_1-a_2)\varphi(XY,Z)+(a_5-a_3)\varphi(YZ,X)\\
+(a_6-a_4)\varphi(ZX,Y)+(a_1-a_4)\varphi(X,YZ)+(a_5-a_2)\varphi(Y,ZX)+(a_6-a_3)\varphi(Z,XY).
\end{array}
\end{equation}

Remark that   we reobtain  Lemma \ref{2} considering $a_1=a_2=a_3=a_4=a_5=a_6$.

%If $\dim A=1$, since $\mu_0 \neq 0$, any bilinear form $\varphi$ satisfies $\varphi=a\mu_0$ with $a \in \K$. and
%$\delta^2_{AA,\mu_0}\varphi \circ \Phi_w=0$ with no assumption on $w$.
%Consider now  that $\dim A \geq 2$. Equation  (\ref{pp}) with $X=Y=Z$ reduces to $$(a_1-a_2-a_3-a_4-a_5+a_6)(\varphi(XX,X)+\varphi(X,XX))=0.$$
 %Since $A^2 \neq 0$, we can choose $X$ such as $XX\neq 0$. If $\varphi$ is such that  $\varphi(XX,X)+\varphi(X,XX)\neq 0$ we obtain
%$$a_1-a_2-a_3-a_4+a_5+a_6=0$$
%and $w=a_1Id+a_2\tau_{12}
%+a_3\tau_{13}+a_4\tau_{23}+a_5c+(a_2+a_3+a_4-a_1-a_5)c^2$.  

%Then from Theorem \ref{co}, $v_{Lad} \notin F_{w}$. 

From the $\mathbb{K}[\Sigma_3]$-module structure of $F_w$ we can deduce \cite{G.R.Nonass} that  
 $v_{3Pa} \in F_{w}$ if only if $F_w$ is odd-dimensional.
In term of the coefficients $a_i$, 
$$\begin{array}{rl}
v_{3Pa} \in F_{w} &  \Leftrightarrow  \exists v=\lambda_1 Id+\lambda_2 \tau_{12}+ \lambda_3  \tau_{13}+\lambda_4  \tau_{23}+\lambda_5 c +\lambda_6 c^2,  v w= v_{3Pa}  \\
&  \Leftrightarrow  M_w V=\ ^t(1,1,1,1,1,1)\\
 \Leftrightarrow &  \left\{
 \begin{array}{l} 
 a_1(\lambda_1-\lambda_5)+a_2(\lambda_2+\lambda_5)+a_3(\lambda_3+\lambda_5)+a_4(\lambda_4+\lambda_5)+a_5(\lambda_6-\lambda_5)=1\\
 a_1(\lambda_2-\lambda_4)+a_2(\lambda_1+\lambda_4)+a_3(\lambda_6+\lambda_4)+a_4(\lambda_5+\lambda_4)+a_5(\lambda_3-\lambda_4)=1\\
 a_1(\lambda_3-\lambda_2)+a_2(\lambda_5+\lambda_2)+a_3(\lambda_1+\lambda_2)+a_4(\lambda_6+\lambda_2)+a_5(\lambda_4-\lambda_2)=1\\
 a_1(\lambda_4-\lambda_3)+a_2(\lambda_6+\lambda_3)+a_3(\lambda_5+\lambda_3)+a_4(\lambda_1+\lambda_3)+a_5(\lambda_2-\lambda_3)=1\\
 a_1(\lambda_5-\lambda_6)+a_2(\lambda_3+\lambda_6)+a_3(\lambda_4+\lambda_6)+a_4(\lambda_2+\lambda_6)+a_5(\lambda_1-\lambda_6)=1\\
 a_1(\lambda_6-\lambda_1)+a_2(\lambda_4+\lambda_1)+a_3(\lambda_2+\lambda_1)+a_4(\lambda_3+\lambda_1)+a_5(\lambda_5-\lambda_1)=1\\
 \end{array}
 \right.
  \end{array}
$$
It implies $(a_2+a_3+a_4)(\lambda_1+\lambda_2+\lambda_3+\lambda_4+\lambda_5+\lambda_6)=3$
so $a_2+a_3+a_4 \neq 0$ and there is a non trivial solution for example $\lambda_5=\lambda_6=\lambda_1=\lambda_2=\lambda_4=\lambda_3=\ds \frac{1}{2(a_2+a_3+a_4)}$, that is $v=\ds \frac{1}{2(a_2+a_3+a_4)}( Id+\tau_{12}+  \tau_{13}+\tau_{23}+c+c^2)$ satisfies $v w=v_{3Pa} $ and  $v_{3Pa} \in F_{w} $. 
%Thus
%For example if $\varphi$ is symmetric with $\varphi(XX,X)\neq 0$ for some $X \in A$, the fact that $\delta^2_{AA,\mu_0}\varphi \circ \Phi_w=0$  implies  $v_{Lad} \notin F_{w}$.
\medskip

\begin{theorem}
Let $(A,\mu_0)$ be a anti-commutative algebra with $\mu_0 \neq 0$ and $dim(A)\geq 2. $
Consider $w=a_1Id+a_2\tau_{12}
+a_3\tau_{13}+a_4\tau_{23}+a_5c+a_6c^2$  a vector of $\KS$ 
We have   that 
 $ \delta^2_{AA,\mu_0}\varphi \circ \Phi_{w}= 0$ for every bilinear map $\varphi$ if
 $a_1-a_2-a_3-a_4+a_5+a_6=0$ and $a_2+a_3+a_4\neq 0$.

\end{theorem}

%%%%%%%%%%%%%%%

\noindent{\bf Remark: deformation quantization and polarization.}  In the following section we recall the notion of polarization-depolarization of a product of an algebra. We will see that when we apply this process to an anti-associative algebra $(A,\mu)$, the associated  skew-symmetric map $\psi$ defined by $\psi(x,y)=\mu(x,y)-\mu(y,x)$ and  the symmetric map $\rho$ defined by $\rho(x,y)=\mu(x,y)+\mu(y,x)$, which provides $(A,\rho)$ with a Jacobi-Jordan algebra structure, are tied up with the graded Leibniz identity. We develop  this point of view in the last section. An algebraic and detailed study of general Jacobi-Jordan  algebras is given in \cite{BF}.  These algebras are also called mock-Lie algebras \cite{Camacho}.

\subsection{Left Leibniz algebras}
Recall that $(A,\mu)$ is a left Leibniz algebra if $\mu$ satisfies the quadratic relation
$$\mu(x,\mu(y, z)) = \mu(\mu(x, y ), z) + \mu(y, \mu(x, z))$$
for any $x,y,z \in A$ what is also written
$$\mathcal{A}_{\mu}(x,y,z)+ \mu(y, \mu(x, z))=0.$$

Let $(A,\mu_0)$ be a commutative left Leibniz algebra. Writing $\mu_0(x,y)=xy$, we have
$$(xy)z-x(yz)+y(xz)=0$$
and $v= Id, w=Id-\tau_{12}.$
For such multiplication we have
$$\delta^2_{v,w,\mu_0}\varphi(x,y,z)=-x\varphi(y,z)+y\varphi(x,z)+z\varphi(x,y)-\varphi(x,yz)+\varphi(y,xz)+\varphi(xy,z)$$
and if $ \delta^2_{v,w,\mu_0}\varphi \circ \Phi_{v_1} =0$ for any bilinear map $\varphi$ then $v_1=0.$
In fact, if $v_1=(a_1,\cdots,a_6)$, then $ \delta^2_{v,w,\mu_0}\varphi \circ \Phi_{v_1} =0$ considered as a linear equation on the formal variables $$x\varphi(y,z),y\varphi(x,z),z\varphi(x,y),\varphi(x,yz),\varphi(y,xz),\varphi(xy,z),$$ each one of the coefficients of these variables being $0$, implies $a_i=0$ for $i=1,\cdots,6$.
So, if $\mu_t$ is a left Leibniz-formal deformation (that is a $(Id,Id-\tau_{12})$-formal deformation) of $\mu_0$, the relation
$$\varphi_1 \bullet^L_{Id} \varphi_1-\varphi_1 \bullet^R_{Id-\tau_{12}} \varphi_1 + \delta^{2,L}_{Id,\mu_0}\varphi_2
- \delta^{2,R}_{Id-\tau_{12},\mu_0}\varphi_2=0$$
can not be reduced.

\begin{remark}\label{rk}
 Maybe there exists some bilinear maps $\varphi_2$ such that $\delta^2_{v,w,\mu_0}\varphi_2 \circ \Phi_{v_1} =0 $ for some $v_1 \neq 0$. But $\varphi_2$ have to satisfy the system associated with the deformation equation $\mu_t \bullet \mu_t =0$ which is very complicated to solve. For these reasons we consider that the equation $\delta^2_{v,w,\mu_0}\varphi_2 \circ \Phi_{v_1} =0 $  is solved for any bililear map $\varphi_2$. 
\end{remark}
\medskip

Let us consider now the relation 
$$\delta^{2,L}_{v,\mu_0}\varphi_1- \delta^{2,R}_{w,\mu_0}\varphi_1 =0$$
which is the order $1$ consequence of the fact that $\mu_t$ is a $(Id,Id-\tau_{12})$ deformation of $\mu_0$.
For any vector $v_1$ we have
$$(\delta^{2,L}_{v,\mu_0}\varphi_1- \delta^{2,R}_{w,\mu_0}\varphi_1)\circ \Phi_{v_1} =0.$$
If $v_1=(a_1,a_2,a_3,-a_1,-a_2,-a_3)$, this equation gives:
$$\begin{array}{rr}
a_1x\psi_1(y,z)+a_2y\psi_1(x,z)-a_3z\psi_1(x,y)+(a_3-a_2)\psi_1(x,yz)&  \\
-(a_1+a_3)\psi_1(y,xz) +(a_1+a_2)\psi_1(z,xy) & =0
\end{array}$$
what is also written 
$$
\begin{array}{l}
a_1(x\psi_1(y,z)-\psi_1(y,xz)+\psi_1(z,xy))\\
+a_2(y\psi_1(x,z)-\psi_1(x,yz)+\psi_1(z,xy))\\
+a_3(-z\psi_1(x,y)+\psi_1(x,yz)-\psi_1(y,xz))= 0               à.
\end{array}
$$
and this is equivalent to the relation
$$x\psi_1(y,z)-\psi_1(y,xz)+\psi_1(z,xy)=0.$$
\begin{proposition}
Let $(A,\mu_0)$ be a commutative left Leibniz algebra.  It is a $(Id,Id -\tau_{12})$-algebra and any $(Id,Id- \tau_{12})$-formal deformation of $\mu_0$ is a deformation quantization of a left-pseudo-Poisson algebra $(A,\mu_0,\psi_1)$ that is a $\K$-vector space $A$ and two bilinear maps satisfying
\begin{enumerate}
\item $\mu_0$ is a commutative left Leibniz multiplication,
\item $\psi_1$ is a skew-symmetric multiplication
\item we have the pseudo Leibniz relation
$$x\psi_1(y,z)-\psi_1(y,xz)-\psi_1(xy,z)=0$$
for any $x,y,z \in A$.
\end{enumerate}
\end{proposition}

\subsection{Right Leibniz algebras}

A right Leibniz algebra $(A,\mu)$, sometimes just called Leibniz algebra, corresponds to the quadratic relation
$$(xy)z-x(yz)-(xz)y=0$$
with $xy$ for $\mu(x,y)$. It is a $(v,w)$ algebra with $v=Id-\tau_{23}$ and $w=Id$. Let us note that if $(u,v) \ra uv$ is a left Leibniz product, then $(u,v) \ra vu$ is a right Leibniz product.

Let us consider a commutative right Leibniz algebra. In this case, the corresponding operator $\delta_{v,w,\mu_0}$ is given by
$$\delta^2_{v,w,\mu_0}\varphi(x,y,z)=z\varphi(x,y)+\varphi(xy,z)-x\varphi(y,z)-\varphi(x,yz)-y\varphi(x,z)-\varphi(xz,y).$$
As before, using the vector $v_1=(a_1,-a_1,a_3,a_4,-a_3,-a_4)$, the equation $\delta^2_{v,w,\mu_0}\varphi(x,y,z) \circ \Phi_{v_1}=0$ is reduced to 
$$z\psi_1(x,y)-\psi_1(x,yz)+\psi_1(y,xz)=0.$$
\begin{proposition}
Let $(A,\mu_0)$ be a commutative Leibniz algebra.  It is a $(Id- \tau_{23},Id)$-algebra and any $(Id- \tau_{23},Id)$-formal deformation of $\mu_0$ is a deformation quantization of a pseudo-Poisson algebra $(A,\mu_0,\psi_1)$ that is
\begin{enumerate}
\item $\mu_0$ be a commutative Leibniz multiplication,
\item $\psi_1$ is a skew-symmetric multiplication
\item we have the pseudo Leibniz relation
$$z\psi_1(x,y)-\psi_1(x,yz)+\psi_1(y,xz)=0$$
for any $x,y,z \in A$.
\end{enumerate}
\end{proposition}

\subsection{Symmetric Leibniz algebras}
A symmetric Leibniz algebra is an algebra $(A, \mu)$ such that for any $x,y,z \in A$, we have
$$
\begin{array}{l}
    (xy)z- x(yz)+y(xz)=0,    \\
   (xy)z-x(yz)-(xz)y=0
\end{array}
$$
with $\mu(x,y)=xy$. Then a symmetric Leibniz algebras is an algebra that is both left Leibniz and right Leibniz. We deduce immediately that if $\mu_t$ is a symmetric Leibniz deformation of a  commutative symmetric Leibniz algebra then 
$$x\psi_1(y,z)-\psi_1(y,xz)-\psi_1(xy,z)=0, \ \psi_1(x,yz)-\psi_1(y,xz)-z\psi_1(x,y)=0$$ where $\psi_1$ is the skew symmetric map associated with $\varphi_1$ and $\mu_0(x,y)=xy$. 
In particular, we have
$$x\psi_1(y,z)+y\psi_1(z,x)+z\psi_1(x,y)=0.$$

\noindent{\bf Remark: Symmetric Leibniz algebras and weakly associative algebras.} 
In \cite{RWass2} we have proved that symmetric Leibniz algebras are weakly associative algebras. 
If $(A, \mu_0)$ is a symmetric Leibniz algebra and if we denote by $\mathcal{A}_{\mu_0}$ the associator of the multiplication $ \mu_0$,  the first identity corresponds to
$$\mathcal{A}_{\mu_0}(x,y,z)=-y(xz)$$
and the second to
$$\mathcal{A}_{\mu_0}(x,y,z)=(xz)y.$$
We deduce that  $(A, \mu_0)$ is a symmetric Leibniz algebra if and only if
$$
\left\{
\begin{array}{l}
\mathcal{A}_{\mu_0}(x,y,z)=-y(xz)\\
 \mathcal{A}_{\mu_0}(x,y,z)=(xz)y.
\end{array}
\right.
$$
In particular we deduce
$$\mathcal{A}_{\mu_0}(x,y,z)+\mathcal{A}_{\mu_0}(y,z,x)=(yx)z-y(xz)
=\mathcal{A}_{\mu_0}(y,x,z)$$
and $(A,\mu_0)$ is also a  weakly associative algebra.
\begin{proposition}
Any symmetric Leibniz algebra is weakly associative.
\end{proposition}
As a consequence, we can consider weakly associative formal deformation of a symmetric Leibniz algebra and, in this case, we find the result of Section 3.5.  
 
 As the symmetric Leibniz are weakly associative, they are Lie-admissible so if  $(A,\mu)$ is symmetric Leibniz the algebra $(A,\psi_\mu)$ is a Lie algebra where $\psi_\mu$ is the skew-symmetric map associated with $\mu$.
 
%Let us note also that $\delta^2_{\mu_0}\varphi \circ \Phi_{v_{Lad}}=0$ if $\mu_0$ is a left or right Leibniz multiplication. 
Consider the vectors $v=Id,w=Id-\tau_{12}, v'=Id-\tau_{23},w'=Id$ and a $(v,w,v',w')$-deformation of $\mu_0$, that is a $(v,w)$-deformation of the $(v,w)$-algebra associated with the equation $(xy)z-x(yz)+y(xz)=0$ and a $(v',w')$-deformation of the $(v',w')$-algebra associated with the equation $(xy)z-x(yz)-(xz)y=0$. Since  $\mu_0$ is a symmetric Leibniz multiplication, the equations coming from the order 2 of a $(v,w,v',w')$-deformation of $\mu_0$ are
$$\delta^2_{v,w,\mu_0,L\ }\varphi_2+\varphi_1 \bullet_{v,w} \varphi_1=0$$
and $$\delta^2_{v',w', \mu_0,R\ }\varphi_2+ \varphi_1 \bullet_{v',w'} \varphi_1=0$$
with 
$$\delta^2_{v,w,\mu_0,L\ }\varphi(x,y,z)=x\varphi(y,z)-y\varphi(x,z)-z\varphi(x,y)+\varphi(x,yz)-\varphi(y,xz)-\varphi(xy,z)$$
and
$$\delta^2_{v',w',\mu_0,R\ }\varphi(x,y,z)=z\varphi(x,y)+\varphi(xy,z)-x\varphi(y,z)-\varphi(x,yz)-y\varphi(x,z)-\varphi(xz,y),$$
The equation $$[\delta^2_{v,w,\mu_0,L\ }\varphi_2+\varphi _1\bullet_{v,w} \varphi_1]\circ \Phi_{Id-\tau_{23}}+[\delta^2_{v',w', \mu_0,R\ }\varphi_2+ \varphi_1 \bullet_{v',w'} \varphi_1] \circ \Phi_{\tau_{13}-c}=0$$
implies $\varphi_1 \bullet_{v_{Lad}} \varphi_1=\varphi_1 \bullet \varphi_1 \circ \Phi_{v_{Lad}}=0.$
We deduce that  if $\mu_0$ is a Leibniz multiplication then $\varphi_1$ is Lie-admissible. As a consequence, if $\mu_0$ is a symmetric Leibniz multiplication, then $\psi_1$ is a Lie bracket.
\begin{proposition}
Let $(A,\mu_0)$ be a commutative symmetric Leibniz algebra.  Any symmetric Leibniz-formal deformation of $\mu_0$ is a deformation quantization of a pseudo-Poisson algebra $(A,\mu_0,\psi_1)$ that is
\begin{enumerate}
\item $\mu_0$ is a commutative Leibniz multiplication,
\item $\psi_1$ is a Lie bracket,
\item we have the pseudo-Leibniz relation
$$\psi_1(x,yz)-\psi_1(y,xz)-z\psi_1(x,y)=0$$
for any $x,y,z \in A$.
\end{enumerate}
\end{proposition}
\medskip

%%%%%%%%%%%%%%%.   Polarization and depolarization.  %%%%%%%%%%%%

\section{Polarization and depolarization of $v$-associative algebras}
Any multiplication $\mu : A 
 \otimes A \to A$ defined by a bilinear application can decomposed into the sum of a 
 commutative multiplication $\rho$ and a skew-symmetric one 
 $\psi$ via the {\em polarization \/} defined by 
 \begin{equation}
\label{pol}
 \rho(x,y) =\frac{1}{2}(x y+y  x)\ 
 \mbox { and }\ 
 \psi(x,y) =\frac{1}{2}(x  y-y  x),\ 
 \mbox { for }\ x,y \in A. 
 \end{equation}
where $\mu(x,y)$ is denoted by $xy$.
 The inverse process of {\em depolarization\/} assembles a commutative 
 multiplication $\rho$ with a skew-symmetric multiplication $\psi$ into the multiplication $\mu$ defined  for any $ x,y \in V$ by 
 \begin{equation}
\label{depol}
 \mu(x,y) =\rho (x, y)+\psi(x,y).
 \end{equation}

 In the following section we first give the well known associative example to illustrate the 
 (de)polarization trick before investigating some other classes of algebras. 
 
\subsection{Associative case \cite{RM}} Assume that $(A,\mu)$ is an associative algebra.
 If we polarize the multiplication $\mu$
 it writes $\mu(x,y)=\rho(x,y)+\psi(x,y)$ and the associativity condition becomes
 equivalent to the following two axioms: 
 \begin{eqnarray} 
 \label{eq:3} 
 \psi(x,\rho(y ,z))&=& \rho(\psi(x,y), z) +\rho(y,\psi(x,z)), 
 \\ 
 \label{eq:1} 
 \psi(y,\psi(x,z)) &=&-\mathcal{A}_\rho (x,y,z). 
 \end{eqnarray}
 To verify this, observe that associativity is equivalent to
$$R(\rho,\psi)=(\psi+\rho)\circ(Id\otimes \psi-\psi\otimes Id)+(\psi+\rho)(Id \otimes  \rho- \rho \otimes Id)=0.$$ 
Moreover it implies that  $R(\rho,\psi) \circ \Phi_v=0$ for any $v \in \KS$. In particular
$$R(\rho,\psi) \circ \Phi_{v_{Lad}}=-4\psi \circ (\psi \otimes Id) \circ \Phi_{Id+c+c^2}= 0$$
and $\psi$ is a Lie bracket. Because of Relation (\ref{eq:3}) the algebra $(A,\rho,\psi)$ is  a nonassociative Poisson algebra in the general case. It is a Poisson algebra if and only if the Lie bracket $\psi$ is $2$-step nilpotent.

Although the associative case is well known, we want to find a systematic method to solve this case which extends  to the other identities that interest us. Let be $\vv \in \KS$ and let us consider the identity $R(\psi,\rho) \circ \Phi_v=0.$ By grouping the terms $\psi(\Id\otimes \psi)$, $\rho(\Id\otimes \rho)$, $\rho(\Id\otimes \psi)$ and finally $\psi(\Id\otimes \rho)$, the coefficients of each of these terms are given by the matricial product
$$N_v \left(
\begin{array}{c}
a_1\\
a_2\\
a_3\\
a_4\\
a_5\\
a_6
\end{array}
\right)$$
where $N_v$ is the transpose of the matrix
$$ \left(
\begin{array}{cccccccccccc}
1&0&1&1&0&-1&1&0&-1&1&0&1\\
0&-1&-1&0&1&-1&0&-1&1&0&1&1\\
-1&0&-1&-1&0&1&1&0&-1&1&0&1\\
-1&-1&0&1&-1&0&-1&1&0&1&1&0\\
1&1&0&-1&1&0&-1&1&0&1&1&0\\
0&1&1&0&-1&1&0&-1&1&0&1&1
\end{array}
\right)
$$
The rank of $N_v$ is $6$. Let us search the vectors of this space associated with minimal relations, that is to say with a maximum of 0 among these components. We obtain the independent vectors
$$\left\{
\begin{array}{l}
(0,0,0,0,0,0,0,1,-1,1,0,0)\\
(0,0,0,0,0,0,-1,0,1,0,1,0)\\
(0,0,0,0,0,0,1,-1,0,0,0,1)\\
\end{array}
\right.
$$
which correspond to the vector $v=(1,-1,1,1,1,-1)$ and $cv,c^2v$ and the relation
$$\psi(x_1,\rho(x_2,x_3))-\rho(x_2,\psi(x_1,x_3))-\rho(x_3,\psi(x_1,x_2))=0.$$
This relation can be written
$\mathcal{L}(\psi,\rho)=0$ where $\psi$ is a Lie bracket and $\rho$ a commutative (nonassociative) multiplication. 
Similarly, we have the three  independent vectors
$$\left\{
\begin{array}{l}
(0,1,0,-1,0,1,0,0,0,0,0,0)\\
(0,0,1,1,-1,0,0,0,0,0,0,0)\\
(1,0,0,0,1,-1,0,0,0,0,0,0)\\
\end{array}
\right.
$$
which correspond to $v=(-1,-1,1,-1,1,1)$ and $cv,c^2v$ and to the relation
$$\psi(x_1,\psi(x_2,x_3))-\rho(x_3,\rho(x_1,x_2))+\rho(\rho(x_3,x_1),x_2)=0$$
that is 
$$\mathcal{A}_\rho+\mathcal{A}_\psi=0.$$ 
In fact $\mathcal{A}_\psi (x_3,x_1,x_2)=-\psi(x_3,\psi(x_1,x_2))+\psi(\psi(x_3,x_1),x_2)=\psi(x_1,\psi(x_2,x_3)).$
This equation implies that $\psi$ is a Lie bracket.
\begin{proposition} Any associative algebra is associated with  polarization/depolarization principle
 to a algebra $(A,\psi,\rho)$ where $\psi$ is a Lie bracket, $\rho$ a commutative multiplication satisfying
\begin{enumerate}
\item $\mathcal{L}(\psi,\rho)=0$
\item $\mathcal{A}_\psi+\mathcal{A}_\rho=0.$
\end{enumerate}
\end{proposition}
In particular, if $\rho$ is associative, then $(A,\psi,\rho)$ is a Poisson algebra with  $2$-step nilpotent Poisson bracket.

\subsection{Lie-admissible case}Recall that a nonassociative algebra $(A,\mu)$ is Lie-admissible if the skew-symmetric bilinear map $\psi$ 
is a Lie bracket. In this case the polarization principle gives no additional relation.

\subsection{Vinberg algebras}
A nonassociative algebra $(A,\mu)$ is a Vinberg algebra if its associator satisfies
$$\mathcal{A}_\mu(x,y,z)-\mathcal{A}_\mu(y,x,z)=0.$$
Since $v_{Lad} \in Span(\mathcal{O}(Id-\tau_{12}))$, such algebras are Lie-admissible. If we polarize the multiplication $\mu$, we obtain
$$
\begin{array}{l}
\psi(x_3,\psi(x_1,x_2)) +\psi(x_1,\rho(x_2,x_3))-\psi(x_,\rho(x_1,x_3))+\rho(x_1,\psi(x_2,x_3))+\rho(x_2,\psi(x_3,x_1))-
\\
-2\rho(x_3,\psi(x_1,x_2))+\rho(x_1,\rho(x_2,x_3))-\rho(x_2,\rho(x_1,x_3))=0.
\end{array}
$$
This relation is also written as the sum
$$ 
\begin{array}{l}
\left\{ (x_1\bu x_3)\bu x_2 - x_1\bu (x_3\bu x_2) - \lbrack x_3, \lbrack x_1, x_2 \rbrack \rbrack \right\}
+
\left\{
[x_1, x_2]\bu x_3 + x_2\bu [x_1, x_3] - [x_1, x_2\bu x_3]\right\} \\
+
\left\{
[x_2, x_1\bu x_3] - x_1\bu [x_2, x_3] - [x_2, x_1]\bu x_3
\right\}
= 0
\end{array}
$$
 of three terms which vanish separately if the multiplication is associative, where $x\bu y =\rho(x,y)$ and  $ [x,y]=\psi(x,y)$.
\subsection{$\tau_{13}$-algebras}
These are the algebras $(A,\mu)$ defined by the quadratic relation
$$\mathcal{A}_{\mu}(x_1,x_2,x_3)-\mathcal{A}_{\mu}(x_3,x_2,x_1)=0.$$This relation is equivalent to
$$\psi(x_1,\psi(x_2,x_3))-\mathcal{A}_{\rho}(x_3,x_1,x_2)=0$$
and $\psi$ is a Lie bracket. This relation is minimal. It can also be written 
$$\mathcal{A}_{\psi} +\mathcal{A}_{\rho}=0$$

\subsection{$(Id+c+c^2)$-algebras} 
In this case we have
$$\mathcal{A}_{\mu}(x_1,x_2,x_3)+\mathcal{A}_{\mu}(x_2,x_3,x_1)+\mathcal{A}_{\mu}(x_3,x_1,x_2)=0.$$
In this case $\psi$ is a Lie bracket and the polarization principe gives after reduction
$$\psi(x_1,\rho(x_2,x_3))+\psi(x_2,\rho(x_3,x_1))+\psi(x_3,\rho(x_1,x_2))=0.$$
It is an identity similar to Equation  (\ref{r2}) obtained by deformation and $(A,\rho,\psi)$ is a non-associative $(Id+c+c^2)$-Poisson algebra. Let us note also that a $(Id+c+c^2)$-algebra is Lie-admissible and $3$-power-associative. 

\subsection{Weakly associative algebras} 
This class of nonassociative algebras has been studied in \cite{RWass} to extend the notion of deformation quantification for associative commutative algebras. Recall that a nonassociative algebra $(A,\mu)$ is weakly associative if we have
$$\mathcal{A}_{\mu}(x_1,x_2,x_3)+\mathcal{A}_{\mu}(x_2,x_3,x_1)-\mathcal{A}_{\mu}(x_2,x_1,x_3)=0.$$
From \cite{RWass}, this identity is equivalent to:
\begin{enumerate}
\item $\psi$ is a Lie bracket,
\item $\rho$ is a commutative multiplication satisfying 
$$\psi(x_1,\rho(x_2,x_3))-\rho(x_2,\psi(x_1,x_3))-\rho(x_3,\psi(x_1,x_2)).$$
\end{enumerate}
In other words, $(\rho,\psi)$ satisfy the Leibniz identity : $\mathcal{L}(\psi,\rho)=0.$
\begin{proposition}
Let $(A,\mu)$ a weakly associative algebra and $(A,\rho,\psi)$ its polarized version. Then 
\begin{enumerate}
\item $(A,\psi)$ is a Lie algebra,
\item $\rho$ is a commutative multiplication, 
\item the multiplications $\rho$ and $\psi$ are tied up by the Leibniz identity
$$\mathcal{L}(\rho,\psi)=0$$
that is $(A,\rho,\psi)$ is a nonassociative Poisson algebra.
\end{enumerate}
\end{proposition}
\noindent

\noindent{\bf Remark.} If we refer to \cite{G.R.Nonass, RWass}, weak associativity corresponds to a point of the family of nonassociative algebras corresponding to the identity
$$\mathcal{C}_\mu(\alpha)(x,y,z)=2\mathcal{A}_\mu (x,y,z)+(1+\alpha)\mathcal{A}_\mu (y,x,z)+\mathcal{A}_\mu (z,y,x)+\mathcal{A}_\mu (y,z,x)+(1-\alpha)\mathcal{A}_\mu (z,x,y) =0$$
with $\alpha =-1/2.$ In fact, considering the vectors $v=Id-\tau_{12}+c$ and $v'=2Id+\frac{1}{2}\tau_{12}+\tau_{13}+c+\frac{3}{2}c^2$, we have from \cite{G.R.Nonass} $\dim F_{v'}=\dim F_v=4$. Since $(Id +\frac{3}{2}\tau_{13}+\tau_{23}+\frac{3}{2}c+c^2)\circ v=v'$ we deduce $F_v=F_{v'}$. 

\medskip

If we consider the vector $v_1=\frac{1}{3}Id-\tau_{12}+\frac{7}{12}\tau_{13}+\frac{1}{4}c^2$, then the polarization of $\mathcal{C}_{\mu}(\alpha) \circ \Phi_{v_1}=0$ gives the relation
$$\mathcal{L}(\rho,\psi)(x_1,x_2,x_3)-\gamma\psi(x_1,\psi(x_2,x_3))-2\psi(x_3,\psi(x_1,x_2))=0$$
with $\gamma=\frac{2}{3}(2\alpha-1)$. Since $v_1$ is  inversible in the algebra $\KS$, this relation is equivalent to $\mathcal{C}_\mu(\alpha)=0.$ 
\begin{proposition}
Let $v=2Id+(1+\alpha)\tau_{12}+\tau_{13}+c+(1-\alpha)c^2$ with $\alpha \neq 1$. Then any $v$-algebra is Lie-admissible and $3$-power-associative. The relation $\mathcal{A}_\mu \circ \Phi_v =0$ is equivalent to
$$\mathcal{L}(\rho,\psi)(x_1,x_2,x_3)-\gamma(\psi(x_1,\psi(x_2,x_3))-2\psi(x_3,\psi(x_1,x_2))$$
with $\gamma=\frac{2}{3}(2\alpha-1)$. In particular, if $\alpha=\frac{1}{2}$, then  $(A,\mu)$ is weakly associative and  we have in this case
$$\mathcal{L}(\rho,\psi)=0$$
that is $(A,\psi,\rho)$ is a nonassociative Poisson algebra.
\end{proposition}

\noindent {\bf Remark.} Polarization of a 3-power associative algebra $(A,\mu).$ The equation corresponding to 3-power associativity 
$$\mathcal{A}_{\mu}(x_1,x_2,x_3)+\mathcal{A}_{\mu}(x_2,x_1,x_3)+\mathcal{A}_{\mu}(x_1,x_3,x_2)+\mathcal{A}_{\mu}(x_3,x_2,x_1)+\mathcal{A}_{\mu}(x_2,x_3,x_1)+\mathcal{A}_{\mu}(x_3,x_1,x_2)\! =\! 0$$
is equivalent to 
$$\psi(x_1,\rho(x_2,x_3))+\psi(x_2,\rho(x_3,x_1))+\psi(x_3,\rho(x_1,x_2))=0.$$
Thus the polarized version of $(A,\mu)$ is   the algebra $(A,\rho,\psi)$ where the commutative multiplication $\rho$ and skew-symmetric multiplication $\psi$ (which is not a Lie bracket) are linked by a   $(Id+c+c^2)$-Leibniz  rule.

\subsection{Poisson algebras}
The depolarization of Poisson algebras has be made in \cite{RM} and \cite{RWass2}. It reinterprets Poisson algebras as structures with one nonassociative product which satisfies
$$3(x(yz)-(xy)z)=-(xz)y-(yz)x+(yx)z+(zx)y.$$

\section{Polarization of $(\K[\ST])^2$-algebras}
\subsection{Polarization of anti-associative algebras}
 Let $(A,\mu)$ be an anti-associative algebra:
\begin{equation}
\label{aa}
\mathcal{AA}_\mu(x,y,z)=\mu(x,\mu(y,z))+\mu(\mu(x,y),z)=0.
\end{equation}
 Let $\psi$ and $\rho$  be the skew-symmetric and symmetric bilinear maps associated to $\mu$ by the polarization principle. To simplify the notations we put $\mu(x,u)=xy, \psi(x,y)=[x,y], \rho(x,y)= x \bu y$.  The identity (\ref{aa}) is equivalent to 
 $$x\bu (y \bu z)+z  \bu (x\bu y)  +x \bu [y,z] +z \bu [x,y]+ [x,y \bu z]-[z,x\bu y]+[x,[y,z]]-[z,[x,y]]=0.$$
 We denote by $\mathcal{AA}(\bu,[,])(x,y,z)=0$ this identity. 
 If $v=Id+\tau_{13}$ then $\mathcal{AA}(\bu,[\, ,\, ])\circ \Phi_v(x,y,z)=0$ is equivalent to 
 $$x\bu (y \bu z) + z \bu (x \bu y)+ [x,[y,z]]-[z,[x,y]]=0$$ which can also been written
 $$\mathcal{AA}_\bullet +\mathcal{AA}_{[\, , \, ]}=0.$$
 It implies that  $$(\mathcal{AA}_\bullet +\mathcal{AA}_{[\, , \, ]})\circ \Phi_{Id+c+c^2}=0$$
 that is 
 $$x\bu (y \bu z) +y \bu (z \bu x)+ z \bu (x \bu y)=0.$$ 
We obtain the following result:
 \begin{proposition}
 Let $(A,\mu)$ be an \ant algebra. If $\rho$ is the symmetric map attached to $\mu$, the algebra $(A,\rho)$ is a Jacobi-Jordan algebra.
 \end{proposition}
 Recall that these algebras have been studied in \cite{BF} where it is proven that these algebras  are commutative nilalgebra of index at most three and conversely. The authors give also classifications for the dimension less or equal to $5$.
 
 Let us note that the Jacobi-Jordan algebras which arise from \ant algebras satisfy also
 $$x \bu (y \bu (z \bu t))=(x \bu y ) \bu (z \bu t)=0$$
 that is all the product of order $4$ are null (see Section 4.2). 
 
 Let us consider now a vector $v$ with 
 $$a_6=a_1-a_2+a_3, \  a_5=a_1+a_3-a_4.$$
 In this case $\mathcal{AA}(\bu,[,])\circ \Phi_v(x,y,z)=0$  is equivalent to 
 $$
 \begin{array}{l}
 \lambda_1(x \bu [y,z]+y \bu [z, x]-[x,y \bu z]+[y,z \bu x])+ \lambda_2 (z \bu [x,y]+[x, y \bu z]-[y,z \bu x]) \\
 + \lambda_3 (y \bu [z,x]-[x,y \bu z]+[z,x \bu y])=0,
 \end{array}$$
 with $\lambda_1=a_1-a_4, \ \lambda_2=a_1-a_2, \ \lambda_3 = a_3-a_2.$
 This identity is equivalent to
  $$\left\{
 \begin{array}{l}
x \bu [y,z]+y \bu [z, x]-[x,y \bu z]+[y,z \bu x]=0, \\
z \bu [x,y]+[x, y \bu z]-[y,z \bu x]=0, \\
y \bu [z,x]-[x,y \bu z]+[z,x \bu y]=0,
 \end{array}
 \right.$$
or to the axiom
$$x \bu [y,z]+[x \bu y, z]+[y,z \bu x]=0.$$
In fact, if $G(x,y,z)=x \bu [y,z]+[x \bu y, z]+[y,z \bu x]=0,$ then the previous system writes
$$\left\{
 \begin{array}{l}
G(x,y,z)+G(y,z,x)=0, \\
G(z,x,y)=0, \\
G(y,z,x)=0.
 \end{array}
 \right.$$
Then $G(x,y,z)+G(y,z,x)+G(z,x,y)=0$ which gives
$$x \bu [y,z]+y \bu [z,x]+z \bu [x,y]=0.$$
\begin{proposition}
Let $(A,\mu)$ be an \ant  algebra. If $\psi,\rho$ are the skew-symmetric and symmetric maps attached to $\mu$, then
\begin{enumerate}
  \item The algebra $(A,\rho)$ is a Jacobi-Jordan algebra.
  \item This algebra $(A,\rho)$  acts as antiderivation on the skew-symmetric algebra $(A,\psi)$.
\end{enumerate} 
  \end{proposition}
  In fact, if $f_x(y)=x \bu y$, then $x \bu [y,z]= -[x \bu y, z]-[y,z \bu x]$ can be written 
  $$f_x[y,z]=-[f_x(y),z]-[y,f_x(z)].$$
  
 \noindent{\bf Remark on the graded Leibniz identity.} The previous identity can be written with the Leibniz identity considering a degree on the operation:
   \begin{definition}
   If 
$\varrho, \eta$ are two multiplication that are symmetric or skew-symmetric, we consider $|\varrho|$ and $|\eta|$ their degree which is $0$ if the operation is skew-symmetric and $1$ if the operation is skewsymmetric. We call graded Leibniz identity on $\eta$ and $\varrho$ 
  \begin{equation}
  \label{graded Leibniz}
  \mathcal{L}_g( \eta,\varrho)(x,y,z)=\varrho(x , \eta(y, z))+(-1)^{|\varrho|}\eta(y ,\varrho (x, z))+(-1)^{|\varrho|} \eta(\varrho (x,y),z).
  \end{equation}
  \end{definition}
  So we obtain with $\psi$ which is skewsymmetric and $\rho$ which is symmetric,
  \begin{enumerate}
\item   $$\mathcal{L}_g( \rho,\psi )(x,y,z)=\psi(x , \rho(y, z))-\rho(y ,\psi (x, z))- \rho(\psi (x,y),z)=0$$
  that is the classical Leibniz identity $$[x,y\bu z]-y\bu[x,z]-[x,y]\bu z$$
  and the skew-symmetric $(A,\psi)$  acts as a derivation on the symmetric algebra $(A,\rho)$
  \item $$\mathcal{L}_g( \psi,\rho)(x,y,z)=\rho(x , \psi(y, z))+\psi(y ,\rho (x, z))+ \psi(\rho (x,y),z)=0$$
  that is $$x\bu [y,z]+[y, x \bu z]+[x\bu y , z]$$
  and the symmetric $(A,\rho)$  acts as an antiderivation on the skew-symmetric algebra $(A,\psi)$
  \item $$\mathcal{L}_g( \psi,\psi)(x,y,z)=\psi(x , \psi(y, z))-\psi(y ,\psi (x, z))- \psi(\psi (x,y),z)$$
  that is $$[x,[y,z]]+[y,[x,z]]+[z,[x,y]]$$
  and the classical Jacobi equation  $[x,[y,z]]+[y,[x,z]]+[z,[x,y]]=0$ writes $\mathcal{L}_g( \psi,\psi)=0$
  \item $$\mathcal{L}_g( \rho,\rho)(x,y,z)=\rho(x , \rho(y, z))+\rho(y ,\rho (x, z))+ \rho(\rho (x,y),z)$$
  that is $$x\bu(y\bu z)+y\bu(z\bu x)+z\bu(x\bu y)$$
  \end{enumerate}
   Then in the graded case some relations obtained previously  become  natural .
  
  \subsection{A remark on the Jacobi-Jordan algebras}
  In the previous section, we have seen that the operad associated with the \ant algebra was not Koszul implying that the cohomology of deformations of these algebras was the cohomology of the minimal model.  It was maybe interesting to look  at this problem for 
  %the commutative \ant algebra that is 
  the  Jacobi-Jordan algebra (this problem is analogous to compare Hochschild and Harrison cohomologies for associative algebras).  We denote by $\mathcal{JJ}ss$ the operad associated to the Jacobi-Jordan algebras. It is clear that
  $$\dim \mathcal{JJ}ss(2)=1, \ \dim \mathcal{JJ}ss(3)=2.$$
  The vector space $\mathcal{JJ}ss(4)$ is generated by the element $a(b(cd))$ and their images by $Id \otimes c$ where $c$ is a cycle in $\Sigma_3$ and $(ab)(cd)$ and their images by $\tau \otimes \tau$ where $\tau$ is the generator of $\Sigma_2$. Then we obtain $15$ generators.  The commutativity and the Jacobi-Jordan condition imply that we have $10$ independent relations. We deduce
  $$\dim \mathcal{JJ}ss(4)=5.$$
  Recall also that if a quadratic operad $\mathcal{P}$ is Koszul, then its Poincar\'e series $g_{P}(t)=-t+ \sum_{k \geq 2} (-1)^k\frac{\dim P(k)}{k!}t^k$ and the Poincar\'e
series of its dual $\mathcal{P}^!$
are tied by the functional equation 
$g_{P}(-g_{P!}(-t))=t. $  Since $g_{\mathcal{JJ}ss}(t)=-t+ \frac{1}{2}t^2-\frac{1}{3}t^3+ \frac{5}{24}t ^4+ \cdots$, the inverse serie is
$$a(t)=-t+\frac{1}{2}t^2-\frac{1}{6}t^3+\cdots$$ On the other hand, any algebra on the dual operad $\mathcal{JJ}ss^!$ is anti-associative and skew-symmetric. 
This can be viewed computing the ideal of relations of this operad. In fact, if $<,>$ denotes the inner product which defines the dual operad of a quadratic operad, we have $<x_1(x_2x_3),(x_1x_2)x_3>=-1$, $<(x_1x_2)x_3,(x_1x_2)x_3>=1$ and $<x_1(x_2x_3),x_1(x_2x_3)>=1$ implying that $<x_1(x_2x_3)+(x_1x_2)x_3,,x_1(x_2x_3)+(x_1x_2x_3)>=0$ and any $\mathcal{JJ}ss^!$-algebra is anti-associative. We deduce, from the anti-associativity, that $\mathcal{JJ}ss^!(4)=0$. Since we have also
$$\dim \mathcal{JJ}ss^!(2)=1, \ \dim \mathcal{JJ}ss^!(3)= 3,$$
we deduce that the generating serie of $\mathcal{JJ}ss^!$ is 
$$g_{\mathcal{JJ}ss^!}(t)=-t+ \frac{1}{2}t^2-\frac{1}{2}t^3.$$

and cannot be a Poincar\'e series of a quadratic operad. Then
\begin{proposition}
The operad $\mathcal{JJ}ss$ of the commutative Jacobi-Jordan algebras is not Koszul. In particular the cohomology of deformations of a Jacobi-Jordan algebra is the cohomology of the minimal model.
\end{proposition}
The determination of the minimal model is similar to those proposed in \cite{M-R-galgalim}.

%%%%%%%%%%%%%
Let us note that the Koszulness of a quadratic operad can be read on associated free algebras. More precisely, a quadratic operad is Koszul if the corresponding free algebras are Koszul algebras. Let us determine the free Jacobi-Jordan algebra. We denote by $JJ(X)$ the (non unitary) free algebra with one generator. Since the in a Jacobi-Jordan algebra we have $x^3=0$, then
$$JJ(X)=\K\{X,X^2\}.$$
Let us denote by $JJ(X,Y)$ the (non unitary) free algebra with two generators. It is a graded algebra $JJ(X,Y)=\oplus _{k \geq 1}JJ_d(X,Y)$ where $JJ_d(X,Y)$ is the subspace of vectors of degree $k$. We have
$$JJ_1(X,Y)=\K\{X,Y\}, \ JJ_2(X,Y)=\K\{X^2, Y^2, XY\}.$$
To compute $JJ_3(X,Y)$ we consider the terms $$X(XX)=(XX)X=X^3, Y^3, X(XY), Y(XX), Y(YX), X(YY).$$  We know that $X^3=Y^3=0.$ We have also
$$2X(XY)+ X^2Y=0, \ 2Y(XY)+XY^2=0.$$
Then
$$JJ_3(X,Y)=\K\{X^2Y, XY^2\}.$$
Let us consider now the terms of degree $4$. Recall that for the Jacobi-Jordan algebras arising from \ant algebras, any terms of degree at least equal to $4$ is $0$. Let us look now the general case. It is clear that
$$X^2X^2=XX^3=Y^2Y^2=YY^3=XY^3=X^3Y=0.$$
For the other terms, that is $$(XY)(XY), X(X(XY)),X(X^2Y),X(XY^2), X(Y(XY)), Y(Y(XY)),Y(XY^2), Y(X^2Y)$$ we have
\begin{enumerate}
\item $X(X(XY))+X(X(YX))+X(Y(XX))=0$ that is $2X(X(XY))=-X(X^2Y)$,
\item $X(X^2Y)+X^2(XY)+ Y(X^2X)=0$ that is $X(X^2Y)=-X^2(XY)$,
\item $X(X(XY))+X((XY)X)+(XY)(XX)=0$ that is $2X(X(XY))=-X^2(XY).$
\end{enumerate}
We deduce
$$2X(X(XY))=-X(X^2Y)=X^2(XY)=-X^2(XY)$$
that is
$$X(X(XY))=X(X^2Y)=X^2(XY)=0.$$
Likewise
$$Y(Y(XY))=Y(Y^2X)=Y^2(XY)=0.$$
In other words, if a term contains a variable of degree $3$, this term vanishes. As for the other terms, we have
\begin{enumerate}
\item $X(XY^2)+X(Y(XY))+X(XY^2)=0$ that is $2X(XY^2)=-X((XY)Y)$,
\item $X(XY^2)+X(Y^2X)+Y^2(XX)=0$, that is $2X(XY^2)=-X^2Y^2$,
\item $(X(XY))Y+( (XY)Y)X+(YX)(XY)=0$ that is (XY)(XY)=-(X(XY))Y-X((XY)Y)
\end{enumerate}
and we deduce
$$X^2Y^2=-2X(XY^2)=-2 Y(X^2Y)=4X((XY)Y)=4(Y((XY)X)=-2(XY)^2.$$
Then
$$JJ_4(X,Y)=\K\{X^2Y^2\}.$$
Let us consider now the terms of degree $5$. Commutativity allows us to consider only the products schematized by
$$\star(\star(\star (\star \star)), \ \ (\star \star)(\star (\star \star)), \ \  (\star \star)(\star (\star \star)).$$
In the first case, the relation $X(X(XY))=X(X^2Y)=X^2(XY)=0, Y(Y(XY))=Y(Y^2X)=Y^2(XY)=0, X((XY)Y)=aX(XY^2)=bX^2Y^2$
show that these products of degree $5$ are reduced to products of type $X(X^2Y^2)$ which are also null. In the second case we have to compute products of type $X((XY)(XY))$ or $X(X^2(XY))$ or $X(Y^2(XY))$ or $X(X^2Y^2)$. In all these cases these products are $0$. In the third case, the Jacobi-Jordan relation  and the computation
$$(X(X^2Y))X=(X^2Y^2)X=0, \ (Y(X^2Y))Y=(X^2Y^2)Y=0$$
shows also that these products are zero.

Then we have
$$JJ_5(X,Y)=0.$$
\begin{proposition} The free Jacobi-Jordan algebra with two generators $JJ(X,Y)$ if of finite dimension and
$$JJ(X,Y)=\K\{X,Y\} \oplus \K\{X^2,Y^2,XY\} \oplus \K\{X^2Y,XY^2\} \oplus \K\{X^2Y^2\}.$$
\end{proposition}

\medskip

Let us  now look at the dual algebra that is \ant skew-symmetric algebra. We denote by $AAS(X_1, \cdots, X_n)$ the free \ant skew-symmetric algebra with $n$-generators. We know that all the product of degree $4$ are zero. Then $AAS(X_1, \cdots, X_n)$ is a subalgebra of $\sum_{1 \leq k \leq 3}AAS^k(X_1, \cdots, X_n)$. Let us determine these algebras.
\begin{enumerate} 
\item $AAS(X)=\K\{X\}$. In fact by antisymmetry $X^2=0$.
\item $AAS(X,Y)=\K\{X,Y,XY\}$. In fact $X^2=Y^2=0$, $X^3=Y^3=0$ and by antiassociativity, $X(XY)=0.$ In this case all the product of degree $3$ are zero.
\item $AAS(X,Y,Z)=\K\{X,Y,Z,XY,XZ,YZ,X(YZ)\}$. 
\item $AAS(X_1, \cdots, X_n)=\K\{X_1, \cdots, X_n\}\oplus \K\{X_iX_j, 1\leq i < j \leq n\} \oplus \K\{X_i(X_lX_k), 1 \leq  i <j <k \leq n\}$.
\end{enumerate}
The Hilbert serie of a graded algebra $A=\oplus A_k$ is $\sum_{k \geq 0}\dim (A_k)t^k$. Then these series are 
\begin{enumerate}
\item for $AAS(X)$: $1+t$
\item for $AAS(X,Y)$: $1+2t+t^2=(1+t)^2$
\item for $AAS(X,Y,Z)$ : $(1+3t+3t^2+t^3=(1+t)^3$ 
\item for $AAS(X_1, \cdots, X_n), n \geq 4$: $1+nt+ \binom {n}{2} t^2+\binom{n}{3} t^3.$
\end{enumerate}
From \cite{Iyudu}, $1+t$, $(1+t)^2$, $(1+t)^3$ are the Hilbert series of Koszul algebras. For the other cases, this problem will be solved latter.

\subsection{Leibniz algebras}
Recall that a Leibniz algebra is a quadratic algebra whose multiplication $\mu(x,y)=xy$ satisfies the identity
\begin{equation}
\label{ll}
x(yz)-(xy)z-y(xz)=0.
\end{equation}
Let $(\rho,\psi)$ be the pair of bilinear maps given by the polarization of $\mu$. As in the previous case, we write $\mu(x,y)=xy, \psi(x,y)=[x,y], \rho(x,y)=x \bu y$. 
Then (\ref{ll}) is equivalent to
$$\begin{array}{l}
R(\rho,\psi)(x,y,z)=x \bu (y \bu z) -(x \bu y) \bu z - y \bu ( x \bu z) +x \bu [y,z]-[x,y]\bu z - y \bu [x,z] \\
+[x, y \bu z]-[x\bu y,z]-[y, x \bu z]+[x,[y,z]]-[[x,y],z]-[y,[x,z]]=0\\
\end{array}
$$
Let $v=(a_1,a_2,a_3,a_4,a_5,a_6)$ be in $\KS$ where the $a_i$ are the components in the canonical basis. In a first time we consider the vector $v=(1,1,0,0,0,0).$ Then $R(\rho,\psi) \circ \Phi_v$ gives in this case the following identity:
$$x \bu (y\bu z)= [x, y \bu z].$$
This implies
$$
 x \bu [y,z]-[x,y]\bu z - y \bu [x,z] 
+2[x, y \bu z]-2[y, x \bu z]+[x,[y,z]]-[[x,y],z]-[y,[x,z]]=0.
$$
Composing this last identity by $\Phi(v)$ with $v=(a_1,a_2,a_3,a_4,a_1-a_2+a_3,a_1-a_2+a_4)$, we obtain
$$x \bu [y,z]+y\bu [z,x]+z \bu [x,y]+3([x,[y,z]]-[[x,y],z]-[y,[x,z]])=0.$$
Then we obtain
$$[x,y\bu z]+[x \bu z,y]-z \bu [x,y]=J_{[,]} (x,y,z)$$
where $J_{[,]} (x,y,z)=[x,[y,z]]-[[x,y],z]-[y,[x,z]]$ is the Jacobi condition for the skew-symmetric map $\psi$.

\begin{proposition}
Let $(A,\mu)$ a Leibniz algebra. It is associated, from the polarization - depolarization principle with a triple $(A,\rho,\psi)$ where $\rho(x,y)= x \bu y$ is a commutative multiplication, on $A$, $\psi(x,y)=[x,y]$ a skew-symmetric multiplication on $A$ satisfying
\begin{enumerate}
 \item $x \bu [y,z]-[x,y]\bu z =[x, y \bu z]$,
\item $[x,y\bu z]+[x \bu z,y]-z \bu [x,y]=J({[,]}) (x,y,z)$
\end{enumerate}
where $J(\psi)$ is the Jacobiator ( $\psi \circ  (Id \otimes \psi ) \circ \Phi_{Id+c+c^2})$, for any $x,y,z \in A$.
\end{proposition}

\begin{corollary}
With the hypothesis of the previous proposition, assume now that $\psi(x,y)=[x,y]$ is a Lie bracket. In this case, for any $z \in A$ the map $f_z (y)=x \bu y$ is a derivation of the Lie algebra $(A,\psi)$.
\end{corollary}
In fact  $[x,y\bu z]+[x \bu z,y]-z \bu [x,y]=J({[,]}) (x,y,z)$ is reduced to
$$[x,f_z(y)]+[f_z(x),y]-f_z([x,y])=0.$$

\medskip

\noindent{\bf Remark: Case of symmetric Leibniz algebras}.  Recall that such algebras correspond to the two identities:
$$
\left\{
\begin{array}{l}
x(yz)-(xy)z-y(xz)=0,\\
(xy)z-x(yz)-(xz)y=0
\end{array}
\right.
$$
This pair of relations is equivalent to
$$\mathcal{A}_{\mu}(x,y,z)-y(xz)=0, \ \mathcal{A}_{\mu}(y,z,x)+(yx)z=0$$
that implies
$$\mathcal{A}_{\mu}(x,y,z)+ \mathcal{A}_{\mu}(y,z,x)- \mathcal{A}_{\mu}(y,x,z)=0.$$
\begin{proposition}
Any symmetric Leibniz algebra is weakly associative. In particular $(A,\rho,\psi)$ is a nonassociative Poisson algebra. 
\end{proposition}
Let us note that the first half of this proposition is the content of Proposition 29.
\section{Deformation quantization and polarization}
Previous studies show that in many cases, if not almost all, there is a close link between the algebras obtained by the formal deformation process and that of polarization. We will summarize this link.
\medskip

{\tiny
\begin{tabular}{|c|c|c| }
\hline
% after \\ : \hline or \cline{col1-col2} \cline{col3-col4} ...
  Type of algebras & Type of algebras appearing  & Type of algebras appearing  \\
  & in formal deformations & by polarization process \\
 \hline
 Associative  &  associative deformation:  &   \\
 & Poisson algebra &  Nonassociative Poisson algebra \\
& & \\
    &  weakly-associative deformation:  & \\
  & Nonassociative Poisson algebra & \\
  \hline
  Lie-admissible & Lie-admissible algebra  &  Lie algebra \\
  \hline
$(Id+c+c^2)$-associative & Nonassociative $(Id+c+c^2)$-Poisson algebra & Nonassociative $(Id+c+c^2)$-Poisson algebra\\
i.e. $G_5$-associative & & \\
\hline
Vinberg algebras & Lie-admissible algebra with & \\
& $(Id-\tau_{12})$-Leibniz condition &
\\
\hline
Weakly-associative & Nonassociative Poisson & Nonassociative Poisson \\
\hline
Anti-associative & Anti-Poisson algebras = & Anti-Poisson algebra = \\
& Jacobi-Jordan algebra & Jacobi-Jordan algebra\\
\hline
Leibniz & Pseudo-Poisson  & Pseudo-Poisson \\
\hline
Symmetric Leibniz &  Pseudo-Poisson  & Nonassociative Poisson  \\
\hline
\end{tabular}
}

\medskip

\noindent {\bf Acknoledgement} The author gratefully thank  the referee for fruitful comments and remarks

\end{document}